\documentclass[a4paper,11pt]{article}
\usepackage[english]{babel}
\usepackage[psamsfonts]{amssymb}
\usepackage{cmmib57}
\usepackage{exscale}
\usepackage{amsmath}
\usepackage{amscd}
\usepackage{latexsym}
\def\1{{\rm 1 }\hskip -0.21truecm 1}
\newcommand{\al}{\alpha}

\newcommand{\del}{\delta}
\newcommand{\Del}{\Delta}
\newcommand{\ep}{\epsilon}

\newcommand{\Lam}{\Lambda}

\newcommand{\gam}{\gamma}
\newcommand{\ffi}{\varphi}
\newcommand{\h}{\mathcal{H}}
\newcommand{\re}{\mathbb{R}}
\newcommand{\red}{\mathbb{R}^{d}}

\newcommand{\e}{\mathbb{E}}
\newcommand{\beq}{\begin{equation}}
\newcommand{\beqn}{\begin{equation*}}
\newcommand{\eeq}{\end{equation}}
\newcommand{\eeqn}{\end{equation*}}
\newcommand{\tm}{[2^mt]}
\newcommand{\tmm}{[2^mt_1]}
\newcommand{\tmmm}{[2^mt_2]}

\newcommand{\qed}{\quad{$\square$}\medskip }
\newtheorem{Th}{Theorem}
\newtheorem{prop}[Th]{\bf Proposition}

\newtheorem{lema}[Th]{ {\bf{Lemma}} }
\newtheorem{corol}[Th]{ {\bf{Corollary}} }
\newtheorem{teor}[Th]{ {\bf{Theorem}} }

\begin{document}
\begin{titlepage}
\mbox{}
\begin{center}
{\LARGE \bf Large deviations for rough paths of the\\[1ex]
fractional Brownian motion}

\vspace*{2cm}
\small
\begin{tabular}{l@{\hspace{2cm}}l}
{\small\sc Annie Millet}$^{\,*,\dagger}$&{\small\sc Marta Sanz-Sol\'e}$^{\,*}$\\
Laboratoire de Probabilit\'es et &Facultat de Matem\`atiques\\
Mod\`eles Al\'eatoires,&Universitat de Barcelona\\
Universit\'es Paris 6 - Paris 7, &Gran Via 585\\
Bo\^{i}te courrier 188, 4, Place Jussieu&E-08007 Barcelona\\ 
F-75252 Paris Cedex 5, and& marta.sanz@ub.edu\\
SAMOS-MATISSE, Universit\'e Paris 1,&\\
amil@ccr.jussieu.fr&
\end{tabular}
\end{center}

\vspace*{3cm}

\noindent{\bf Abstract:} Starting from the construction of a geometric rough path associated
with a fractional Brownian motion with Hurst parameter $H\in]\frac{1}{4},\frac{1}{2}[$ given
by Coutin and Qian in \cite{coutinquian}, we prove a large deviation principle in the space of geometric
rough paths, extending classical results on Gaussian processes. As a by-product, geometric rough
paths associated to elements of the reproducing kernel Hilbert space of the fractional Brownian motion
are obtained and an explicit integral representation is given. 



\mbox{}

\vfill

\begin{itemize}
\item[{$^*$}]{\small Supported by the grants BFM 2003-01345, HF2003-006 from the Direcci\'on
General de Investigaci\'on, Ministerio de Educaci\'on y Ciencia, Spain.}\\[-.8cm]
\item[{$^\dagger$}]{\small Supported by the program SAB 2003-0082 from the
Direcci\'on General de
Universidades, Ministerio de Educaci\'on y Ciencia, Spain.}
\end{itemize}
\end{titlepage}
\section{\bf Introduction}

In the seminal paper \cite{lyons}, Lyons developed a sophisticated mathematical theory 
to analyse dynamical systems with an external rough force acting
as a control and influencing their evolution. One of the key ideas 
is to keep the non-commutative structure of controls on small time steps.
Rough controls are constructed as elements of direct sums of tensor spaces endowed with
a topology associated with the $p$-variation distance. Dynamical systems are proved to be continuous 
functionals of their rough path controls with respect to this topology. This result is called
the {\it universal limit theorem}.

Stochastic modeling deals basically with rough path controls. Indeed, the
ground-breaking It\^o's theory on stochastic differential equations is based on
Brownian motion, which has almost surely nowhere differentiable sample paths but only
$\alpha$-H\"older continuous ones, with $\alpha\in]0,\frac{1}{2}[$. 
Note that the solution of a multidimensional stochastic
It\^o's differential equation is not a continuous functional of the driving Brownian motion.
From Lyons perspective, the
rough path character of Brownian motion is caught by increments of both, its trajectories 
and those of the L\'evy area process. 
His approach provides a kind of pathwise calculus well-suited for system control in a
stochastic context. We refer the reader to \cite {lq} and \cite{l}, where the basic ingredients
of the theory are presented.

It\^o's theory has been extensively developed in many different directions,
including finite and infinite dimensional settings. Recently, increasing
attention is being devoted to a particular stochastic control rougher than
the Brownian motion: the fractional Brownian motion with Hurst parameter
$H\in ]0,\frac{1}{2}[$. Unlike the classical Brownian process ($H=\frac{1}{2}$), the fractional Brownian
motion does not have independent increments and possesses long-range memory. Many
problems in traffic networks, hydrology and economics, just to mention a few examples, share these
properties and therefore can be  realistically analysed including this process in their mathematical
formulation. In \cite{m}  a large survey on fractional Brownian motion is given. Some of the recent
developments  concerning fractional Brownian motion are
employed in this paper (see for instance \cite{amn}, \cite{coutinquian}, \cite{ccm}, \cite{decreuse}).
These references contain an exhaustive list of contributors to the subject and are suggested to those
who would like to have a broad picture on the subject.
\smallskip

In this article, we are interested in the rough
path associated with a fractional Brownian motion with Hurst parameter $H\in]\frac{1}{4},\frac{1}{2}[$,
constructed in \cite{coutinquian}. The main goal has been to establish a large deviation
principle.  For $H=\frac{1}{2}$, this question
has been addressed in \cite{lqz} and the possibility of the extension given in our work is mentioned.
However, we believe that it is not a straightforward one and gives rise to interesting
mathematical issues which need new ideas to be solved satisfactory.
For values of $H$ in $]\frac{1}{2},1[$, the problem has an almost obvious answer -see the remark following
the proof of Proposition \ref{p2}.  

In order to give a more detailed description of the results
in their context, some basic notions on rough paths analysis and some notation should
be set up.

Let $T>0$ and $\mathbb{B}$ be a Banach space. For $p\ge 1$, the $p$-variation norm of a
function $x:[0,T]\longrightarrow\mathbb{B}$ is defined by
\beqn
||x||_p=\left(\sup_{\mathcal {P}}\sum_l|x_{t_l}-x_{t_{l-1}}|^p_{\mathbb{B}}\right)^{\frac{1}{p}},
\eeqn
where the supremum runs over all finite partitions $\mathcal{P}$ of $[0,T]$.
In the sequel we shall take $T=1$ and consider $\mathbb{B}=\red$. 

A continuous map $X$ defined on the simplex $\Del=\{(s,t): 0\le s\le t\le 1\}$, taking values on
the truncated tensor algebra
\beqn
T^{[p]}(\red)= \re\oplus\red\oplus(\red)^{\otimes 2}\oplus\cdots\oplus(\red)^{\otimes [p]}
\eeqn
is called a {\it rough path} in $T^{[p]}(\red)$ of roughness $p$, if
$X_{s,t}=(1,X^1_{s,t},\cdots,X^{[p]}_{s,t})$, $(s,t)\in\Del$, satisfies the properties:
\begin{description}
\item {(a)} Finite $p$-variation: $\max_{1\le j\le [p]}\left(\sup_{\mathcal{P}}\sum_{l}
|X^j_{t_{l-1},t_l}|^{\frac{p}{j}}\right)^{\frac{j}{p}}<\infty$.
\item{(b)} Multiplicative property: $X_{s,t}=X_{s,u}\otimes X_{u,t}$, for any $(s,u), (u,t)\in\Del$.
\end{description}
The set of rough paths in $T^{[p]}(\red)$ is a metric space with the $p$-variation distance
\beq
\label{pvarm}
d_p(X,Y)= \max_{1\le j\le [p]}\left(\sup_{\mathcal{P}}\sum_{l}
|X^j_{t_{l-1},t_l}-Y^j_{t_{l-1},t_l}|^{\frac{p}{j}}\right)^{\frac{j}{p}}.
\eeq
Assume that the function $x$ has finite total variation. For any $j=1,\dots,[p]$, $(s,t)\in\Del$,
consider the $j$-th iterated integral
\beq
\label{mint}
X_{s,t}^j = \int\cdots\int_{s<t_1<\cdots<t_j<t} dx_{t_1}\otimes\cdots\otimes dx_{t_j}.
\eeq
It is easy to check that $X_{s,t}=(1,X^1_{s,t},\cdots,X^{[p]}_{s,t})$ defined in this way is a rough
path. We shall refer to this class of objects as {\it smooth rough paths} lying above $x$.

The space of {\it geometric rough paths} with roughness $p$ is the closure
of the set of smooth rough paths with respect to the $p$-variation metric. An important class
in stochastic analysis of geometric rough paths are those obtained from smooth rough paths based
on linear interpolations of $x$. They shall be denoted by $\mathcal{D}_p(\red)$. Indeed, linear
interpolations of interesting examples like Brownian motion, $\mathbb{B}$-valued Wiener process, free
Brownian motion and fractional Brownian motion have been successfully used to define the corresponding
geometric rough path (see \cite{lq}, \cite{llq}, \cite{cd}, \cite{coutinquian}, respectively).
\medskip

In this paper, we consider a $d$-dimensional fractional Brownian motion $W^H = (W^H_t, t\in[0,1])$
with Hurst parameter $H\in]0,1[$. Its reproducing kernel
Hilbert space, denoted by $\h^H$, consists of functions $h: [0,1]\rightarrow \red$ that can be
represented as
\beq
\label{h}
h(t) = \int_0^t K^H(t,s) \dot h(s) ds,
\eeq
where $K^H(t,s)$ is the kernel defined by
\begin{align}
K^H(t,s)=&\, c_H(t-s)^{H-\frac{1}{2}}\nonumber\\
&+c_H\left(\frac{1}{2}-H\right) \int_s^t(u-s)^{H-\frac{3}{2}}
\left(1-\left(\frac{s}{u}\right)^{\frac{1}{2}-H}\right) du,
\label{1.0}
\end{align}
for $c_H>0$,  $0< s < t \leq 1$, and $\dot h\in L^2([0,1])$.
The scalar product in $\h^H$ is given by
\beqn
\langle h_1,h_2\rangle_{\h^H}=\langle\dot h_1,\dot h_2\rangle_{L^2[0,1]}
\eeqn
(see \cite{decreuse}, Theorem 3.3).

For $0 < s<t$ 
\beq
\label{a.1}
\frac{\partial K^H}{\partial t}(t,s) = c_H\, \left(H-\frac{1}{2}\right)\,
\left(\frac{s}{t}\right)^{\frac{1}{2}-H}(t-s)^{H-\frac{3}{2}}.
\eeq
Note that for $H\in]0,\frac{1}{2}[$,  
$|K^H|(dt,s) = -\frac{\partial K^H}{\partial t}(t,s)\,
1_{]s,1[}(t)\, dt$.
\smallskip

Let $\e= \mathcal{C}_0([0,1]; \red)$, endowed with the topology of the supremum norm and let
$P^H$ be the law of $W^H$ on $\e$. The triple $(\e,\h^H,P^H)$ is an abstract Wiener space.
We shall denote by $i^H$ the continuous dense embedding of $\h^H$ into $\e$.

A classical result of the theory of Gaussian processes
(see for instance \cite{deuschel}, Theorem 3.4.12)  establishes that the family $(\ep P^H, \ep>0)$
of Gaussian probabilities satisfies a large deviation principle on $\e$ with good rate function
\beq
\label{1.1}
{\Lam}^H (x)=\begin{cases} \frac{1}{2} ||(i^H)^{-1}(x)||^2_{\h^H}& \text{if $x\in i^H(\h^H)$},\\
+\infty&\text{otherwise}.
\end{cases}
\eeq

Along this article  we deal with values of $H$ in $]\frac{1}{4},\frac{1}{2}[$.
For the sake of simplicity, we shall skip any reference to the parameter $H$ in the sequel
and write $W$ instead of $W^H$, $\h$ instead of $\h^H$, etc.

For any $m\in\mathbb{N}$, we consider the $m$-th dyadic grid $(t_l^m=l2^{-m}, l=0,1,\dots,2^m)$ and
set $\Delta^m_lW = W_{t_l^m}-W_{t_{l-1}^m}$, for any $l=1, \dots,2^m$.

Denote by
$W(m)=(W(m)_t, t\in[0,1])$ the process obtained by linear interpolation of $W$ on the $m$-th dyadic grid.
That is,
$W(m)_0=0$ and for $t\in]t_{l-1}^m,t_l^m]$,
\beq
\label{1.2}
W(m)_t=W_{t_{l-1}^m}+2^m(t-t_{l-1}^m)\ \Delta^m_lW.
\eeq
Let $p\in]1,4[$ be such that $Hp>1$. In \cite{coutinquian}, a geometric rough path
with roughness $p$, lying above $W$
is obtained as a limit in the $p$-variation distance (\ref{pvarm}) of the sequence of smooth rough paths
$F(W(m))=(1,W(m)^1,W(m)^2,W(m)^3)$ defined as in (\ref{mint}).
We denote this object by $F(W)$. By its very construction, $F(W)\in \mathcal{D}_p(\red)$.

As has been mentioned before, our purpose is to establish a large deviation principle on
$\mathcal{D}_p(\red)$ for the family of probability laws of $(F(\ep W), \ep\in(0,1))$, extending the
classical Schilder result for Gaussian processes.
By means of the universal limit theorem of \cite{lyons}, the result can be transferred to 
stochastic differential equations driven by fractional Brownian motion.

The next section 2 is devoted to the proof of the main result. We follow the same strategy as in \cite{lqz}.
That is, since the smooth rough paths based on linear interpolations of the process $W$ are easily seen to
satisfy a large deviation principle, we only need to prove that they are {\it exponentially good}
approximations of $W$. In comparison with \cite{lqz}, there are essentially two new difficulties coming
up. Firstly, time increments of fractional Brownian motion are not independent and secondly, we need
to deal with third order geometric rough paths, making arguments a bit more involved.
The main tools to be used are the hypercontractivity inequality for Gaussian chaos (see \cite{lt}) and
a collection of covariance type estimates for $W$ proved in
\cite{coutinquian}. As a by-product, we prove the existence of a geometric rough path associated
with each element $h$ in the reproducing kernel Hilbert space $\h$. Section three is entirely
devoted to give a precise description of this geometric rough path in terms of indefinite multiple
integrals. The results might be understood as  deterministic versions of those given in \cite{amn}
for stochastic integrals with respect to Gaussian Volterra processes (see also \cite{decreuse2}).
In our case, integrands and integrators are of Volterra type, because of the representation (\ref{h}).
The interest of these results goes beyond the framework of this work; they shall be useful in
the characterization of the topological support of the law of the rough path associated with the fractional
Brownian motion.

As is being usual, we denote throughout the proofs different constants by the same letter.

\section{\bf The large deviation principle}

We want to prove the following.
\begin{teor}
\label{t1}
Let $H\in]\frac{1}{4},\frac{1}{2}[$, $p\in]1,4[$ be such that $Hp>1$. The family
of probability laws of $(F(\ep W), \ep\in(0,1))$ satisfies a large deviation principle on
$\mathcal{D}_p(\red)$ with the good rate function defined for $X\in \mathcal{D}_p(\red)$ by:
\beq
\label{2.1}
\mathcal{I}(X)=\begin{cases} \frac{1}{2} ||i^{-1}(X^1_{0,.})||^2_{\h}& \text{if $X^1_{0,.}
\in i(\h)$},
\\ +\infty&\text{otherwise}.
\end{cases}
\eeq
\end{teor}

Let us start by setting the method of the proof, that we borrow from \cite{lqz} and fix the notations
to be used in the sequel.

Let $Z(m)=(W_{t_l^m}, 1\le l\le 2^m)$. Clearly, $Z(m)=\Phi_m\circ W$, with $\Phi_m: \e\rightarrow
(\re^d)^{2^m}$ a continuous map. The explicit form of the smooth rough path lying above $W(m)$ shows that
there exists a continuous map $\Psi_m:(\re^d)^{2^m}\rightarrow \mathcal{D}_p(\red)$ such that
$F(W(m))= (\Psi_m\circ\Phi_m)(W)$. Consequently, the contraction principle implies that for any $m$ the family
of probability laws of $\left(F(\ep W(m)), \ep\in(0,1)\right)$ satisfies a large deviation principle on
$\mathcal{D}_p(\red)$ with the good rate function
\beq
\label{2.2}
\mathcal{I}_m(X)= \inf\{\Lam (x): x\in \e, (\Psi_m\circ\Phi_m)(x)= X\},
\eeq
$X\in \mathcal{D}_p(\red)$.

We then transfer the large deviation principle from $F(\ep W(m))$
to $F(\ep W)$. At first we shall prove that for any $\del>0$, \beq
\label{2.3} \lim_{m\to\infty}\limsup_{\ep\to 0}\ep^2\log
P\left(d_p\left(F(\ep W(m)),F(\ep W)\right)>\del\right) = -\infty.
\eeq

For any $h\in\h$, let $h(m)$ denote the smooth function obtained by linear interpolation of $h$ on
the $m$-th dyadic grid and let $F(h(m))$ be the corresponding smooth rough path.

We will prove that for every $\al>0$,
\beq
\label{2.4}
\lim_{m,m'\to\infty}\sup_{||h||_{\h}\le \al} d_p\left(F(h(m)),F(h(m'))\right) =0.
\eeq
This result gives in particular the existence of a geometric rough path $F(h)$
in
$\mathcal{D}_p(\red)$ obtained as the limit in the $d_p$-variation distance of $F(h(m))$.
In the last part of the article we shall identify $F(h)$ as a triple of integrals.

By means of an extension of the contraction principle (see \cite{dz}, Theorem 4.2.23),
(\ref{2.3}), (\ref{2.4}) provide a proof of Theorem \ref{t1}.

Let us introduce some technicalities to deal with the $p$-variation distance $d_p$.

If $X,Y$ are rough paths of degree $[p]$, we set for $j=1,\dots, [p]$,  $\gam>0$,
\beqn
D_{j,p}(X,Y) = \left(\sum_{n=1}^\infty n^\gam
\sum_{l=1}^{2^n}|X_{t^n_{l-1},t^n_l}^j-Y_{t^n_{l-1},t^n_l}^j|^{\frac{p}{j}}\right)^{\frac{j}{p}}
\eeqn
and $D_{j,p}(X)=D_{j,p}(X,0)$.

Owing to results proved in \cite{llq} and \cite{coutinquian} (see also \cite{lq}), for any
$p\in ]3,4[$, $\gam>p-1$,
\begin{align}
d_p(X,Y)&\le C\,\max\Big( D_{1,p}(X,Y), D_{1,p}(X,Y)\big[ D_{1,p}(X)+D_{1,p}(Y)\big],
D_{2,p}(X,Y),\nonumber\\ 
&D_{2,p}(X,Y)\big[D_{1,p}(X)+D_{1,p}(Y)\big],
D_{1,p}(X,Y)\big[D_{2,p}(X)+D_{2,p}(Y)\big],\nonumber\\
&D_{1,p}(X,Y)\big[D_{1,p}(X)^2+D_{1,p}(Y)^2\Big], D_{3,p}(X,Y)\Big). \label{2.5}
\end{align}

Therefore, similar
arguments as in  \cite{lqz}, pp. 273-274 show that (\ref{2.3}) follows
from the following statement.

\begin{prop}
\label{p1}
Let $p\in]1,4[$ be such that $pH>1$. Then,

\noindent (a) For any $j=1,2,3$, there exists a sequence $c_j(m)$ converging to zero as $m$ tends to infinity
such that for every $q>p$,
\beq
\label{2.6}
\big(E\left(D_{j,p}(W(m),W)^q\right)\big)^{\frac{1}{q}} \le c_j(m) q^{\frac{j}{2}}.
\eeq
(b)  For any $j=1,2$, there exists a constant $c_j$ such that for every $q>p$,
\beq
\label{2.7}
\sup_{m\in\mathbb{N}}\left(E\left(D_{j,p}(W(m))^q\right)\right)^{\frac{1}{q}}\le c_jq^{\frac{j}{2}}.
\eeq
\end{prop}

\noindent{\bf Proof:} We shall denote by $g$ a standard normal random variable and observe that,
as a consequence of the hypercontractivity inequality (see e.g. \cite{lt}, page 65),
$(E|g|^q)^{\frac{1}{q}}\le (q-1)^{\frac{1}{2}}$, for any $q\in]2,\infty[$.
Along the proof, for any $n\geq m$, $l=1,\dots,2^n$, we denote by
$k:=k(n,m,l)$ the unique integer in $\{1,2,\dots,2^m\}$ such that
\beq
\label{k}
t_{k-1}^m\le t^n_{l-1}<t^n_l<t_k^m .
\eeq

\noindent{\it First order terms.} Let $j=1$.  From the definition of $W(m)^1$ it follows easily,
\beqn
D_{1,p}(W(m),W) = \Big(\sum_{n=m+1}^\infty
n^{\gam}\sum_{l=1}^{2^n}|2^{m-n}\Delta_k^mW-\Delta^n_lW|^p\Big)^{\frac{1}{p}}.
\eeqn
As in \cite{lqz}, for $q>p$ and $m\geq 0$ set $A(m,q)=\left(\sum_{n=m+1}^\infty
2^n\left(\frac{n^\gam}{a_n}\right)^{\frac{q}{q-p}}\right)^{\frac{q-p}{p}}$, for some
sequence of real positive numbers $(a_n, n\ge 1)$ to be
chosen later.

H\"older's inequality yields
\begin{align}
E\big(D_{1,p}(&W(m),W)^q\big)\le A(m,q)\sum_{n=m+1}^\infty a_n^{\frac{q}{p}}
\sum_{l=1}^{2^n}E\big(|2^{m-n}\Delta_k^mW-\Delta^n_lW|^q\big)\nonumber\\
&\le A(m,q)(2d)^q q^{\frac{q}{2}}\sum_{n=m+1}^\infty
 a_n^{\frac{q}{p}} 2^n  \big(2^{-nq}2^{mq(1-H)}+2^{-nqH}\big)\nonumber\\
&\le A(m,q)(2d)^q q^{\frac{q}{2}}\sum_{n=m+1}^\infty a_n^{\frac{q}{p}}2^{n(1-qH)}.\label{2.8}
\end{align}
Set $a_n= 2^{np(H-\frac{1}{q}-\ep)}$ for some $\ep\in]0,\frac{1}{2}(H-\frac{1}{p})[$;
then the series $ \sum_n a_n^{\frac{q}{p}}2^{n(1-Hq)} $ converges.
Fix $\del>0$
such that $n^\gam\le c2^{n\del p}$ for some $c>0$
and $\ep+\del\in]0,\frac{1}{2}(H-\frac{1}{p})[$.
Then,
\beqn
A(m,q)^{\frac{1}{q}} \le c^{\frac{1}{p}} 2^{-m(H-\frac{1}{p}-\ep-\del)}.
\eeqn
Consequently, $\sup_{q>p} A(m,q)^{\frac{1}{q}}$ tends to zero as $m$ tends to infinity.
By virtue of (\ref{2.8}), the upper bound (\ref{2.6}) for $j=1$ holds true.

The proof of (\ref{2.7}) for $j=1$ is similar. Indeed, proceeding as for the proof of
(\ref{2.8}), we have
\beq
\label{2.8.1}
E\big(D_{1,p}(W)^q\big)\le A(0,q)(2d)^q q^{\frac{q}{2}}\sum_{n=1}^\infty
a_n^{\frac{q}{p}}2^{n(1-qH)}.
\eeq
Since  $\sup_{q>p}A(0,q) < \infty$, the inequalities (\ref{2.8}), (\ref{2.8.1}) yield
(\ref{2.7}).
\smallskip

\noindent{\it Second order terms.} Let $j=2$. For $l=1,\dots,2^n$ set
\beqn
T_2(n,m,l)= W(m+1)_{t^n_{l-1}, t^n_l}^2-W(m)_{t^n_{l-1}, t^n_l}^2.
\eeqn
Assume first $n< m$. 
Quoting equation (20) in \cite{coutinquian},
\beqn
T_2(n,m,l)=\frac{1}{2}\sum_{r=2^{m-n}(l-1)+1}^{2^{m-n}l}
\big(\Delta_{2r-1}^{m+1}W\otimes \Delta_{2r}^{m+1}W
-\Delta_{2r}^{m+1}W\otimes\Delta_{2r-1}^{m+1}W\big).
\eeqn
Clearly, $T_2(n,m,l)=0$ for $d=1$ and for any $d\ge 2$, all the diagonal components $T_2(n,m,l)^{i,i}$ vanish.
Hence, we may assume that $d\ge 2$ and consider only $(i,j)$ components with $i\ne j$. Under these
premises, any couple of random variables $\Delta_k^{m+1}W^i$, $\Delta_h^{m+1}W^j$ are independent.

Owing to the hypercontractivity inequality,
\beqn
\big(E\left|T_2(n,m,l)^{i,j}\right|^q\big)^{\frac{1}{q}}\le C q
\big(E\left|T_2(n,m,l)^{i,j}\right|^2\big)^{\frac{1}{2}}.
\eeqn
Clearly,
\beqn
E\left|T_2(n,m,l)^{i,j}\right|^2 \leq C  \left( T_{21}(n,m,l)^{i,j}+T_{22}(n,m,l)^{i,j}\right),
\eeqn
where
\begin{align}
T_{21}(n,m,l)^{i,j}&= \sum_{r=2^{m-n}(l-1)+1}^{2^{m-n}l}
E\big|\Delta_{2r-1}^{m+1}W_i\; \Delta_{2r}^{m+1}W_j -
\Delta_{2r-1}^{m+1}W_j\; \Delta_{2r}^{m+1}W_i \big|^2 \nonumber\\
&\leq C \sum_{r=2^{m-n}(l-1)+1}^{2^{m-n}l} E\big|\Delta_{2r-1}^{m+1}W_i\big|^2
E\big|\Delta_{2r}^{m+1}W_j\big|^2\nonumber\\
&\le C2^{-n}2^{-4m(H-\frac{1}{4})}.\label{2.90}
\end{align}
Lemma 12 in \cite{coutinquian} yields
\begin{align}
T_{22}(n,m,l)^{i,j}&\le C \sum_{r=2^{m-n}(l-1)+1}^{2^{m-n}l} \left(\sum_{\bar r=r+1}^\infty
(\bar r-r)^{4H-5}\right) 2^{-4(m+1)H}\nonumber\\
&\le C2^{-n}2^{-4m(H-\frac{1}{4})}.\label{2.91}
\end{align}
Consequently,
\beq
\left(E\left|T_2(n,m,l)^{i,j}\right|^q\right)^{\frac{1}{q}}
\le C q 2^{-\frac{n}{2}} 2^{-2m(H-\frac{1}{4})}.
\label{2.9}
\eeq
This inequality holds also true for $n=m$. Indeed, using for instance the identities (5) and
(6) in \cite{lqz} for $m=n+1$ and $m=n$, respectively, we obtain
\beqn
T_2(n,n,l)=\frac{1}{2}\left(\Del_{2l-1}^{n+1}W\otimes\Del_{2l}^{n+1}W-
\Del_{2l}^{n+1}W\otimes\Del_{2l-1}^{n+1}W\right),
\eeqn 
and therefore,
\beqn
\left(E\left|T_2(n,n,l)^{i,j}\right|^q\right)^{\frac{1}{q}}
\le C q 2^{-2nH}.
\eeqn

Fix $M>m$. The above inequality (\ref{2.9}) and Minkowski's inequality imply
\begin{align}
\Big(E\left|W(M)_{t^n_{l-1}, t^n_l}^2-W(m)_{t^n_{l-1},
t^n_l}^2\right|^q\Big)^{\frac{1}{q}} &\le C q
2^{-\frac{n}{2}}\sum_{N=m}^{M-1}2^{-2N(H-\frac{1}{4})}\nonumber\\ &\le C q
2^{-\frac{n}{2}}2^{-2m(H-\frac{1}{4})},\label{2.10}
\end{align}
where $C$ is a constant depending only on $H, p$ and $d$.

By the construction of the rough path lying above $W$, a.s.,
\beqn
\lim_{M\to\infty}W(M)_{t^n_{l-1}, t^n_l}^2 = W_{t^n_{l-1}, t^n_l}^2.
\eeqn
Therefore, Fatou's lemma and (\ref{2.10}) yield for $n\leq m$
\beq
\label{2.11}
\Big(E\left|W_{t^n_{l-1}, t^n_l}^2-W(m)_{t^n_{l-1}, t^n_l}^2\right|^q\Big)^{\frac{1}{q}}
\le C q 2^{-\frac{n}{2}}2^{-2m(H-\frac{1}{4})}.
\eeq

Let $m\le n$;  in this case,
\beqn
W(m)_{t^n_{l-1}, t^n_l}^2 = 2^{2(m-n)-1}(\Delta_k^m W)^{\otimes 2},
\eeqn
where $k=k(n,m,l)$ satisfies (\ref{k}) (see \cite{coutinquian}, equation (17)).
By the hypercontractivity property,
\beq
\label{2.12}
\Big(E\left|W(m)_{t^n_{l-1}, t^n_l}^2\right|^q\Big)^{\frac{1}{q}}
\le C q 2^{-2n} 2^{-2m(H-1)}.
\eeq
The previous estimate (\ref{2.11}) for $n=m$ together with Minkowski's
inequality, imply
\beq
\label{2.13}
\Big(E\left|W_{t^n_{l-1}, t^n_l}^2\right|^q\Big)^{\frac{1}{q}}
\le C q 2^{-2nH}.
\eeq
With (\ref{2.12}) and (\ref{2.13}) we obtain for $m\leq n$,
\beq
\label{2.14}
\Big(E\left|W_{t^n_{l-1}, t^n_l}^2-W(m)_{t^n_{l-1}, t^n_l}^2\right|^q\Big)^{\frac{1}{q}}
\le C q 2^{-2nH}.
\eeq

We now proceed in an analogue way as for $j=1$. For $q>\frac{p}{2}$, set $A_2(q)=\left(\sum_{n=1}^\infty
2^n\left(\frac{n^\gam}{a_n}\right)^{\frac{2q}{2q-p}}\right)^{\frac{2q-p}{p}}$, for some
positive real sequence $(a_n, n\ge 1)$.
By H\"older's inequality,
\beq
\label{2.15}
E\left(D_{2,p}(W(m),W)^q\right)\le A_2(q)\sum_{n=1}^\infty a_n^{\frac{2q}{p}}
\sum_{l=1}^{2^n}E\left|W(m)_{t^n_{l-1}, t^n_l}^2-W_{t^n_{l-1}, t^n_l}^2\right|^q.
\eeq
From (\ref{2.11}) and (\ref{2.14}), it follows that
\begin{align*}
E\Big(D_{2,p}&(W(m),W)^q\Big)\le C A_2(q)q^q\Big[\sum_{n=1}^m a_n^{\frac{2q}{p}}
2^{-n(\frac{q}{2}-1)} 2^{-2mq(H-\frac{1}{4})}\\
&+\sum_{n=m+1}^\infty a_n^{\frac{2q}{p}} 2^{-n(2qH-1)}\Big].
\end{align*}
Notice that, for any $\eta\in]0,2H-\frac{1}{2}[$,
\beq
\label{2.16} \sum_{n=1}^m a_n^{\frac{2q}{p}}2^{-n(\frac{q}{2}-1)}
2^{-2mq(H-\frac{1}{4})} \le 2^{-mq\eta}\sum_{n=1}^\infty
a_n^{\frac{2q}{p}}2^{-n[q(2H-\eta)-1]}. \eeq

Let $a_n=
2^{-np(\ep-H+\frac{\eta}{2}+\frac{1}{2q})}$, with $\ep>0$. Then
the series 
$ \sum_n  a_n^{\frac{2q}{p}}2^{-n[q(2H-\eta)-1]}$ converges. Moreover, this
choice of $a_n$ yields
\begin{align*}
&A_2(q)=\left( \sum_{n=1}^\infty n^{\frac{2\gam q}{2q-p}}
2^{-\frac{npq}{2q-p}(-2\ep+2H-\eta-\frac{2}{p})}\right)^{\frac{2q-p}{p}},\\
&\sum_{n=m+1}^\infty a_n^{\frac{2q}{p}} 2^{-n(2qH-1)} = \sum_{n=m+1}^\infty 2^{-nq(\eta+2\ep)}.
\end{align*}
Let $\eta, \ep  $ and $\delta$ be positive reals such that
 $\delta+\ep+\frac{\eta}{2}<H-\frac{1}{p}$, and $n^\gam\le C2^{np\delta}$, for
some
$C>0$. Then $\sup_{q>\frac{p}{2}}\left(A_2(q)\right)^{\frac{1}{q}} <\infty$ and consequently,
\beq
\label{2.17}
\big(E\left(D_{2,p}(W(m),W)^q\right)\big)^{\frac{1}{q}} \le C q 2^{-m\eta},
\eeq
proving (\ref{2.6}) for $j=2$.

By a similar approach, using the estimate (\ref{2.13}), we can prove that
\beqn
\big(E\left(D_{2,p}(W)^q\right)\big)^{\frac{1}{q}}\le C q.
\eeqn
Thus, (\ref{2.7}) for $j=2$ holds true.
\smallskip

\noindent{\it Third order terms.} Finally, let us prove (\ref{2.6}) for
$j=3$. Assume first $n\le m$; then for any $l=1,\dots,2^n$,
\beq
\label{2.18}
E\left|W(m+1)^3_{t_{l-1}^n,t_l^n}-W(m)^3_{t_{l-1}^n,t_l^n}\right|^2
\le C 2^{-n(1+2H)} 2^{-m(4H-1)}. \eeq
Indeed, for $n<m$, the
inequality is proved in \cite{coutinquian}, pg. 128. Let $n=m$;
quoting \cite{coutinquian} pg. 119, for any $n\ge 1$, we write
\begin{align}
W(n+1)^3_{t_{l-1}^n,t_l^n} = & \sum_{k=2l-1}^{2l} \Big( W(n+1)^3_{t_{k-1}^{n+1},t_k^{n+1}}\nonumber\\
&+ W(n+1)^1_{t_{l-1}^n,t_{k-1}^{n+1}}\otimes W(n+1)^2_{t_{k-1}^{n+1},t_k^{n+1}}\nonumber\\
&+W(n+1)^2_{t_{l-1}^{n},t_{k-1}^{n+1}}\otimes W(n+1)^1_{t_{k-1}^{n+1},t_{k}^{n+1}}\Big).
\label{2.18.1}
\end{align}
Fix $k\in\{2l-1, 2l\}$; it is easy to check that for any $q\in[2,\infty)$,
\beq
\label{2.19}
\Big(E\left|W(n+1)^1_{t_{l-1}^{n},t_{k-1}^{n+1}}\right|^q\Big)^{\frac{1}{q}}
+\Big(E\left|W(n+1)^1_{t_{k-1}^{n+1},t_k^{n+1}}\right|^q\Big)^{\frac{1}{q}}
\le C q^{\frac{1}{2}} 2^{-(n+1)H}.
\eeq
Applying (\ref{2.12}) we obtain
\beq
\label{2.20}
\Big(E\left|W(n+1)^2_{t_{k-1}^{n+1},t_k^{n+1}}\right|^q\Big)^{\frac{1}{q}}
+\Big(E\left|W(n+1)^2_{t_{l-1}^{n},t_{k-1}^{n+1}}\right|^q\Big)^{\frac{1}{q}}
\le C q 2^{-(2n+2)H}.
\eeq
Moreover, for any $m\le n$,
\beq
\label{2.21}
W(m)^3_{t_{l-1}^n,t_l^n}= \frac{2^{3(m-n)}}{3!}(\Delta_{k}^m W)^{\otimes 3},
\eeq
with $k = k(n,m,l)$ satisfying (\ref{k}).
Since $W(m)^3_{t_{l-1}^n,t_l^n}$ belongs to the third order Gaussian chaos,
the hypercontractivity property yields for $m\leq n$,
 \beq
\label{2.22}
\left(E\left|W(m)^3_{t_{l-1}^n,t_l^n}\right|^q\right)^{\frac{1}{q}}
\le C q^{\frac{3}{2}}2^{-3n}2^{-3m(H-1)}.
\eeq
From (\ref{2.18.1}) - (\ref{2.20}) and (\ref{2.22}), we obtain
\beq
\label{2.22.1}
E\left|W(n+1)^3_{t_{l-1}^n,t_l^n}\right|^2 \le C 2^{-6nH}.
\eeq
This upper bound, together with (\ref{2.22}) for $m=n$ and $q=2$, imply
the validity of (\ref{2.18}) for $n=m$.

By virtue of the hypercontractivity property and (\ref{2.18}) we deduce for $n\leq m$,
\beqn
\left(E\left|W(m+1)^3_{t_{l-1}^n,t_l^n}-W(m)^3_{t_{l-1}^n,t_l^n}\right|^q\right)^{\frac{1}{q}}
\le C q^{\frac{3}{2}}2^{-n(\frac{1}{2}+H)} 2^{-m(2H-\frac{1}{2})}.
\eeqn
Hence, Minkowski's inequality yields
\beqn
\left(E\left|W(M)^3_{t_{l-1}^n,t_l^n}-W(m)^3_{t_{l-1}^n,t_l^n}\right|^q\right)^{\frac{1}{q}}
\le C q^{\frac{3}{2}}2^{-n(\frac{1}{2}+H)} 2^{-m(2H-\frac{1}{2})}
\eeqn
for any $M>m\geq n$.

We observe that,  a.s. $\lim_{M\to\infty}W(M)^3_{t_{l-1}^n,t_l^n}=W^3_{t_{l-1}^n,t_l^n}$. Therefore,
Fatou's Lemma yields for $m\geq n$,
\beq
\left(E\left|W^3_{t_{l-1}^n,t_l^n}-W(m)^3_{t_{l-1}^n,t_l^n}\right|^q\right)^{\frac{1}{q}}
\le C q^{\frac{3}{2}}2^{-n(\frac{1}{2}+H)} 2^{-m(2H-\frac{1}{2})}.
\label{2.23}
\eeq

Suppose $m\le n$. Applying the previous estimate (\ref{2.23}) and (\ref{2.22}) with $m=n$,
we obtain
\beqn
\left(E\left|W^3_{t_{l-1}^n,t_l^n}\right|^q\right)^{\frac{1}{q}}
\le C q^{\frac{3}{2}} 2^{-3nH}.
\eeqn
Therefore, using again (\ref{2.22}) we deduce for $m\leq n$,
\beq
\label{2.24}
\left(E\left|W^3_{t_{l-1}^n,t_l^n}-W(m)^3_{t_{l-1}^n,t_l^n}\right|^q\right)^{\frac{1}{q}}
\le C q^{\frac{3}{2}}2^{-3nH}.
\eeq

For $q>\frac{p}{3}$, let $A_3(q)= \left(\sum_{n=1}^\infty
2^n\left(\frac{n^\gam}{a_n}\right)^{\frac{3q}{3q-p}}\right)^{\frac{3q-p}{p}}$, where $a_n$ is a
sequence of positive numbers to be determined later. 
H\"older's inequality yields
\beqn
E\big(D_{3,p}\left(W(m),W\right)^q\big) \le A_3(q) \sum_{n=1}^\infty a_n^{\frac{3q}{p}}
\sum_{l=1}^{2^n}E\left|W(m)^3_{t_{l-1}^{n},t_{l}^{n}}-W^3_{t_{l-1}^{n},t_{l}^{n}}\right|^q.
\eeqn
By means of (\ref{2.23}), (\ref{2.24}) we obtain,
\begin{align*}
E\big(D_{3,p} & \left(W(m),W\right)^q\big) \le A_3(q) q^{\frac{3q}{2}} \Big(\sum_{n=1}^m a_n^{\frac{3q}{p}}
2^{-nq(\frac{1}{2}+H-\frac{1}{q})-mq(2H-\frac{1}{2})}\\
&\quad \quad +\sum_{n=m+1}^\infty a_n^{\frac{3q}{p}} 2^{-n(3qH-1)} \Big).
\end{align*}
Let $\eta\in]0,2H-\frac{1}{2}[$; clearly,
\beqn
\sum_{n=1}^m a_n^{\frac{3q}{p}}2^{-nq(\frac{1}{2}+H-\frac{1}{q})-mq(2H-\frac{1}{2})}
\le 2^{-mq\eta}\sum_{n=1}^m a_n^{\frac{3q}{p}}2^{-nq(3H-\frac{1}{q}-\eta)}.
\eeqn
Set $a_n= 2^{-np(\frac{\ep}{3}-H+\frac{\eta}{3}+\frac{1}{3q})}$, with $\ep>0$. Then the series
$\sum_{n=1}^\infty a_n^{\frac{3q}{p}}2^{-nq(3H-\frac{1}{q}-\eta)}$ converges.
Furthermore,
\begin{align*}
&A_3(q) = \left(\sum_{n=1}^\infty n^{\frac{3\gam q}{3q-p}}
2^{-\frac{npq}{3q-p}(-\ep+3H-\eta-\frac{3}{p})}\right)^{\frac{3q-p}{p}}\, ,\\
&\sum_{n=m+1}^\infty
a_n^{\frac{3q}{p}} 2^{-n(3qH-1)}=\sum_{n=m+1}^\infty 2^{-nq(\ep+\eta)}\, .
\end{align*}
Let $\eta>0, \ep >0$ and $ \delta>0$ be such that $3 \del+\ep+\eta<3H-\frac{3}{p}$
and $n^\gam<C2^{np\del}$, for some
$C>0$. Then $\sup_{q>\frac{p}{3}} \left(A_3(q)\right)^{\frac{1}{q}}<\infty$. Thus,
\beqn
\big(E\left(D_{3p}(W(m),W)^q\right)\big)^{\frac{1}{q}} \le Cq 2^{-m\eta},
\eeqn
proving (\ref{2.6}) for $j=3$.
This concludes  the proof of the Proposition. \hfill \qed
\bigskip

In the sequel, we make the convention $K(t,s)=0$ if $s\ge t$, and therefore
write
\beqn
h(t)= \int_0^1 K(t,s) \dot h(s) ds,
\eeqn
for any $h\in\h$.
We denote by $||\cdot||_2$ the usual Hilbert norm in $L^2([0,1])$.

\begin{lema}
\label{l1}
Let $h\in\h$ and $t,t'\in[0,1]$. Then
\beq \label{2.25.pre} |h(t) - h(t')|\le ||\dot h||_2 \; |t-t'|^H.
\eeq In particular, for any $\alpha >0$, \beq \label{2.25}
\sup_{||h||_{\h}\le \alpha}|h(t) - h(t')|\le \alpha |t-t'|^H. \eeq
\end{lema}
\noindent{\bf Proof:} With the above convention on the kernel $K$ and by virtue of Schwarz's inequality,
we have
\begin{eqnarray*}
|h(t) - h(t')|^2&=&\left|\int_0^1 \left(K(t,s)-K(t',s)\right) \dot h(s) ds\right|^2\\
&\leq &||\dot h||_2^2 \int_0^1 \left(K(t,s)-K(t',s)\right)^2 ds\\
&=&||\dot h||_2^2 \; E\left|W_t-W_{t'}\right|^2
=||\dot h||_2^2 \; |t-t'|^{2H}.\qquad \square 
\end{eqnarray*}
\medskip

In the remaining part of the section, $h$ shall denote a fixed element in $\h$
 and $h(m)$, $m\ge 1$, the function
obtained by linear interpolation of $h$ on the $m$-th dyadic grid $(t_l^m=l2^{-m}, l=0,1,\dots,2^m)$.
That is, $h(m)_0=0$ and for $t\in]t_{l-1}^m,t_l^m]$,
\beq
\label{h(m)}
h(m)_t=h(t_{l-1}^m)+2^m(t-t_{l-1}^m)\Delta_l^mh.
\eeq
We shall quote several times algebraic identities set up in \cite{coutinquian} for the processes $W(m)$,
$m\ge 1$, and replace $W(m)$ by $h(m)$. Indeed, their proof rely only on the structure of
the linear interpolations and not on the probabilistic properties of the fractional Brownian motion.
\smallskip

Our next purpose is to prove the convergence stated in (\ref{2.4}). By the inequality (\ref{2.5}),
this amounts to prove the next proposition.
\begin{prop}
\label{p2}
Let $p\in]1,4[$ be such that $pH>1$ and $\alpha >0$. Then,

\noindent (a) For  every $j=1,2,3$,
\beq
\label{2.26}
\lim_{m,m'\to\infty} \sup_{||h||_{\h}\le \alpha}D_{j,p}\left(h(m),h(m')\right)=0.
\eeq
(b) For every $j=1,2$,
\beq
\label{2.27}
\sup_{m\in\mathbb{N}}\sup_{||h||_{\h}\le \alpha} D_{j,p}\left(h(m)\right)<\infty.
\eeq
\end{prop}

\noindent{\bf Proof:}
{\it First order terms.} Let $j=1$ and $k$ the  index satisfying (\ref{k}). By Lemma \ref{l1}, 
\begin{eqnarray*}
\sup_{||h||_{\h}\le \alpha}\left(D_{1,p}\left(h(m),h\right)^p\right)&=&
\sup_{||h||_{\h}\le \alpha}\sum_{n=m+1}^\infty n^\gam \sum_{l=1}^{2^n}|2^{m-n}\Del_k^mh-\Del_l^nh|^p\\
&\le& C\alpha^p \sum_{n=m+1}^\infty n^\gam \sum_{l=1}^{2^n} \big(2^{-mp(H-1)-np}+2^{-npH}\big)\\
&\le& C \alpha^p 2^{-m(pH-1-\ep)},
\end{eqnarray*}
with $\ep\in]0,pH-1[$.
Hence, (\ref{2.26}) holds for $j=1$.

Similarly, for any $\ep\in]0,pH-1[$,
\begin{align*}
\sup_{||h||_{\h}\le \alpha}\left(D_{1,p}(h)^p\right)&=
\sup_{||h||_{\h}\le \alpha}\sum_{n=1}^\infty n^\gam \sum_{l=1}^{2^n}|\Del_l^n h|^p\\
&\le \alpha^p \sum_{n=1}^\infty 2^{-n(pH-1-\ep)} \le C\alpha^p,
\end{align*}
which together with (\ref{2.26}) for $j=1$ give (\ref{2.27}) for $j=1$.
\smallskip

\noindent{\it Second order terms.} Consider now the case $j=2$. Assume first $m \le n$.
 Following \cite{coutinquian}, equation (17),
pg. 118 for $w(m):=h(m)$, and using Lemma \ref{l1}, we have for $m<n$,
\begin{align}
\sup_{||h||_{\h} \le \alpha} &\left|h(m+1)^2_{t_{l-1}^n,t_l^n}-h(m)^2_{t_{l-1}^n,t_l^n}\right|\nonumber \\
&\le \sup_{||h||_{\h}\le
\alpha}\left(\left|h(m+1)^2_{t_{l-1}^n,t_l^n}\right|+\left|h(m)^2_{t_{l-1}^n,t_l^n}\right|\right)\nonumber\\
&\le C2^{2(m-n)}\left(\left|(\Del_k^{m+1}h)^{\otimes 2}\right|+\left|(\Del_k^{m}h)^{\otimes
2}\right|\right)
\nonumber\\
&\le C\alpha^2 2^{-2nH}.
\label{2.28}
\end{align}
Notice that we have also proved that for every $m \le n$,
\beq
\label{2.28.1}
\sup_{||h||_{\h}\le \alpha}\left|h(m)^2_{t_{l-1}^n,t_l^n}\right|\le C\alpha^2 2^{-2nH}.
\eeq
From equation (19) in \cite{coutinquian} pg. 118 and Lemma \ref{l1}, we easily obtain
\beqn
\sup_{||h||_{\h}\le \alpha}\left|h(n+1)^2_{t_{l-1}^n,t_l^n}\right|\le C\alpha^2 2^{-2nH}.
\eeqn
Thus, the above upper bound (\ref{2.28}) holds for any $m \le n$.

Suppose now $n<m$. Quoting \cite{coutinquian}, equation (20), pg
118, we write
\begin{align*}
h(m+1)_{t_{l-1}^n,t_l^n}^2 & -h(m)_{t_{l-1}^n,t_l^n}^2= \frac{1}{2}\sum_{k=2^{m-n}(l-1)+1}^{2^{m-n}l} 
\big(\Del_{2k-1}^{m+1}h\otimes \Del_{2k}^{m+1}h\\
&  - \Del_{2k}^{m+1}h\otimes \Del_{2k-1}^{m+1}h\big),
\end{align*}
for any $l=1,\dots,2^n$.

Fix $d\ge 2$ and components $(i,j)$ of the tensor products with $i\neq j$. Clearly,
\beqn
\left|h(m+1)_{t_{l-1}^n,t_l^n}^{2,i,j}-h(m)_{t_{l-1}^n,t_l^n}^{2,i,j}\right|\le C \left(T^{i,j}_{m,n,l}
+T^{j,i}_{m,n,l}\right),
\eeqn
with
\begin{eqnarray}
T^{i,j}_{m,n,l}& = &
\Big|\sum_{k=2^{m-n}(l-1)+1}^{2^{m-n}l}\Del_{2k-1}^{m+1}h^i\Del_{2k}^{m+1}h^j\Big|\nonumber\\
&=&\Big|\int_0^1\int_0^1\sum_{k=2^{m-n}(l-1)+1}^{2^{m-n}l}\left(K(t_{2k-1}^{m+1},t)
-K(t_{2k-2}^{m+1},t)\right)\nonumber\\
&&\quad \times\left(K(t_{2k}^{m+1},s)-K(t_{2k-1}^{m+1},s)\right)\dot h^i(s)\dot h^j(t) ds dt\Big|.
\label{t}
\end{eqnarray}
Schwarz's inequality yields
\beqn
T^{i,j}_{m,n,l}\le C \alpha^2 \left(T^{i,j}_{m,n,l}(1)+T^{i,j}_{m,n,l}(2)\right)^{\frac{1}{2}},
\eeqn
where
\begin{eqnarray*}
T^{i,j}_{m,n,l}(1)&=&\sum_{k=2^{m-n}(l-1)+1}^{2^{m-n}l}\int_0^1\int_0^1
\left(K(t_{2k-1}^{m+1},t)
-K(t_{2k-2}^{m+1},t)\right)^2\\
&&\quad \times\left(K(t_{2k}^{m+1},s)-K(t_{2k-1}^{m+1},s)\right)^2 ds dt,\\
T^{i,j}_{m,n,l}(2)&=&\sum_{k,k'=2^{m-n}(l-1)+1\atop k<k'}^{2^{m-n}l}\int_0^1\int_0^1
\left(K(t_{2k-1}^{m+1},t)-K(t_{2k-2}^{m+1},t)\right)\\
&&\quad \times \left(K(t_{2k'-1}^{m+1},t)-K(t_{2k'-2}^{m+1},t)\right)
\left(K(t_{2k}^{m+1},s)-K(t_{2k-1}^{m+1},s)\right)\\
&&\quad \times 
\left(K(t_{2k'}^{m+1},s)-K(t_{2k'-1}^{m+1},s)\right) ds dt.
\end{eqnarray*}
Clearly,
\begin{eqnarray}
T^{i,j}_{m,n,l}(1)&=&\sum_{k=2^{m-n}(l-1)+1}^{2^{m-n}l}
E\left|\Delta_{2k-1}^{m+1}W_i\right|^2
 E\left|\Delta_{2k}^{m+1}W_j\right|^2\nonumber\\
&\le & C2^{-m(4H-1)-n}.
\label{2.29}
\end{eqnarray}
Moreover,
\[T^{i,j}_{m,n,l}(2)
=\sum_{k,k'=2^{m-n}(l-1)+1\atop k<k'}^{2^{m-n}l}
E\left(\Del_{2k-1}^{m+1}W_i \; \Del_{2k'-1}^{m+1}W_i\right)\; 
 E\left(\Del_{2k}^{m+1}W_j \; \Del_{2k'}^{m+1}W_j\right).
\]
Notice that, whenever $k<k'$, $2(k'-k)\ge 2$. Hence, applying the upper
bound set up in \cite{coutinquian}, equation (29), pg. 121, we obtain
\begin{eqnarray}
T^{i,j}_{m,n,l}(2)&\le& C 2^{-4(m+1)H}\sum_{k,k'=2^{m-n}(l-1)+1\atop k<k'}^{2^{m-n}l}
|k'-k|^{4(H-1)}\nonumber\\
&\le &C 2^{-4(m+1)H} 2^{m-n}\sum_{r=1}^{2^{m-n}} r^{-4(1-H)}\nonumber\\
&\le & C 2^{-m(4H-1)-n}.
\label{2.30}
\end{eqnarray}
From (\ref{2.29}), (\ref{2.30}), it follows that, if $n<m$,
\beq
\label{2.31}
\sup_{||h||_{\h}\le \alpha}\left|h(m+1)^2_{t_{l-1}^n,t_l^n}-h(m)^2_{t_{l-1}^n,t_l^n}\right|
\le C \alpha^2 2^{-m(2H-\frac{1}{2})-\frac{n}{2}}.
\eeq

Putting together (\ref{2.28}) (valid for $m\le n$) and (\ref{2.31}), we obtain
\begin{align*}
\sup_{||h||_{\h}\le \alpha} \big(D_{2,p} & \left(h(m+1),h(m)\right)\big)^{\frac{p}{2}} \le C\alpha^p
\Big(\sum_{n=1}^{m-1}n^\gam 2^{n(1-\frac{p}{4})}2^{-mp(H-\frac{1}{4})}\\
&+\sum_{n=m}^\infty n^\gam 2^{-n(pH-1)}\Big).
\end{align*}

Let $\ep\in]0,pH-1[$. The above estimates show the existence of some positive real
number $\beta$ such that
\beq
\label{2.32}
\sup_{||h||_{\h}\le \alpha} D_{2,p}\left(h(m+1),h(m)\right)\le C \alpha^2 2^{-m\beta}.
\eeq
This yields (\ref{2.26}) for $j=2$.
The proof of (\ref{2.27}) for $j=2$ is an easy consequence of (\ref{2.32}) and (\ref{2.28.1}).
\smallskip

\noindent{\it Third order terms.}  We finally prove (\ref{2.26}) for $j=3$; note that
these terms only appear when $H\in ]\frac{1}{4},\frac{1}{3}]$, so that $p\in ]3,4[$.

Assume first $m\le n$. In this case, for any
$l=1,\dots,2^n$ and $k$ satisfying (\ref{k})
\beq
\label{2.33}
h(m)^3_{t_{l-1}^n,t_l^n}=\frac{2^{3(m-n)}}{3!}\left(\Delta_k^mh\right)^{\otimes 3}.
\eeq
We shall check that
\beq
\label{2.34}
\sup_{||h||_{\h}\le \alpha}\left|h(m+1)^3_{t_{l-1}^n,t_l^n}-h(m)^3_{t_{l-1}^n,t_l^n}\right|\le C\alpha^3
2^{-3nH}.
\eeq
Indeed, if $m<n$, owing to (\ref{2.33}) and Lemma \ref{l1},
\begin{align*}
\sup_{||h||_{\h}\le \alpha} & \left|h(m+1)^3_{t_{l-1}^n,t_l^n}-h(m)^3_{t_{l-1}^n,t_l^n}\right| \le
\sup_{||h||_{\h}\le \alpha}\Big(\left|h(m+1)^3_{t_{l-1}^n,t_l^n}\right| \\
 &\qquad +\left|h(m)^3_{t_{l-1}^n,t_l^n}\right|\Big)\\
&\le C\alpha^3 2^{-3m(H-1)-3n}\le C\alpha^3 2^{-3nH}.
\end{align*}
For $m=n$, we write the analogue of (\ref{2.18.1}) with $W$ replaced by $h$.  With Lemma \ref{l1},
we can check that each term of the resulting formula is bounded above by $C\alpha^3 2^{-(n+1)3H}$.
Consequently, (\ref{2.34}) holds for any $m\le n$.
\smallskip

Let us now assume $n<m$. We write the identity given in Lemma 11 in \cite{coutinquian} with
$w$ replaced by $h$.
More precisely,
\beqn
\left|h(m+1)^3_{t_{l-1}^n,t_l^n}-h(m)^3_{t_{l-1}^n,t_l^n}\right|\le C\sum_{r=1}^5
\left|I_r(m,n,l)\right|,
\eeqn
with
\begin{align*}
I_1(m,n,l)=&\sum_k\left(h(t_{2k-2}^{m+1})-h(t_{l-1}^n)\right)\otimes \left(\Del_{2k-1}^{m+1}h
\otimes \Del_{2k}^{m+1}h-\Del_{2k}^{m+1}h\otimes\Del_{2k-1}^{m+1}h\right),\\
I_2(m,n,l)=&\sum_{k}\left(\Del_{2k-1}^{m+1}h\otimes\Del_{2k}^{m+1}h-
\Del_{2k}^{m+1}h\otimes \Del_{2k-1}^{m+1}h\right)\otimes \left(h(t_l^n)-h(t_{2k+2}^{m+1})\right),\\
I_3(m,n,l)=&\sum_{k}\Del_{2k-1}^{m+1}h\otimes\left(\Del_{2k}^{m+1}h\otimes \Del_{2k}^{m+1}h
+\Del_{2k-1}^{m+1}h\otimes \Del_{2k}^{m+1}h\right),\\
I_4(m,n,l)=&\sum_{k}\Del_{2k}^{m+1}h\otimes\left(\Del_{2k}^{m+1}h\otimes \Del_{2k-1}^{m+1}h
+\Del_{2k-1}^{m+1}h\otimes \Del_{2k}^{m+1}h\right),\\
I_5(m,n,l)=&\sum_{k}\left(\Del_{2k-1}^{m+1}h\otimes \Del_{2k}^{m+1}h+
\Del_{2k}^{m+1}h\otimes\Del_{2k-1}^{m+1}h\right)\otimes\Del_{2k-1}^{m+1}h,
\end{align*}
where the index $k$ in the sums runs in the set $\{2^{m-n}(l-1)+1,\dots,2^{m-n}l\}$.
The first two terms above have the same structure; the last three ones are also similar.
They shall be analysed separately.

We start with $I_1(m,n,l)$. Notice that if $d=1$, this term vanishes. Moreover, for $d\ge 2$
only the components $I_1(m,n,l)^{\kappa,i,j}$ with $i\ne j$ might not vanish.

Let $i\ne j$. Clearly,
\begin{align*}
\sup_{||h||_{\h}\le \alpha} & \left|I_1(m,n,l)^{\kappa,i,j}\right|\le \sup_{||h||_{\h}\le \alpha}\sup_k 
\left|h(t_{2k-2}^{m+1})-h(t_{l-1}^n)\right|\\
& \quad\quad\times\left(T_{m.n,l}^{i,j}+T_{m.n,l}^{j,i}\right),
\end{align*}
with $T_{m.n,l}^{i,j}$ defined in (\ref{t}). Then, Lemma \ref{l1} together with
(\ref{2.29}) and (\ref{2.30}) yield
\beq
\label{2.35}
\sup_{||h||_{\h}\le \alpha}\left|I_1(m,n,l)^{\kappa,i,j}\right|\le C \alpha^3
2^{-m(2H-\frac{1}{2})-n(H+\frac{1}{2})},
\eeq
and the same estimate holds for $\sup_{||h||_{\h}\le \alpha}\left|I_2(m,n,l)^{\kappa,i,j}\right|$.

Set \beqn
J(m,n,l)=\sum_{k=2^{m-n}(l-1)+1}^{2^{m-n}l}\Del_{a(k)}^{m+1}h\otimes\Del_{b(k)}^{m+1}h
\otimes \Del_{c(k)}^{m+1}h, 
\eeqn
where $a(k), b(k)$ and
$c(k)$ belong to $ \in\{2k-1,2k\}$ and are such that two out of
the three indices agree. Lemma \ref{l1} yields \beqn
\sup_{||h||_{\h}\le \alpha}  |J(m,n,l)| \le C\alpha^3 2^{-m(3H-1)-n},
\eeqn which implies
\beq 
\label{2.36} \sup_{||h||_{\h}\le\alpha}
\left|I_{\mu}(m,n,l)\right|\le C\alpha^3 2^{-m(3H-1)-n}, \eeq
for any
$\mu=3,4,5$.

Since for $n<m$,  $2^{-m(2H+\frac{1}{2})-n(H-\frac{1}{2})}
<2^{-m(3H-1)-n}$, the estimates (\ref{2.35}) and (\ref{2.36}) imply
\beq
\label{2.37}
\sup_{||h||_{\h}\le \alpha}\left|h(m+1)^3_{t_{l-1}^n}-h(m)^3_{t_{l-1}^n}\right|
\le C\alpha^3 2^{-m(3H-1)-n}.
\eeq

By the very definition of $D_{3,p}(h(m+1),h(m))$ and taking into account the results obtained
for $m\le n$ in (\ref{2.34})  and for $m>n$ in (\ref{2.37}),  we obtain
\[ D_{3,p}\big(h(m+1),h(m)\big)^{\frac{p}{3}}\le C \alpha^p\Big(\sum_{n=1}^{m-1}n^\gam
2^{-mp(H-\frac{1}{3})-n\frac{p}{3}} 
+\sum_{n=m}^\infty n^\gam2^{-n(pH-1)}\Big).
\]
Since $p>3$, this yields
\beqn
\sup_{||h||_{\h}\le \alpha}D_{3,p}(h(m+1),h(m)) \le C \alpha^3 2^{-m\beta},
\eeqn
for some real $\beta>0$. This suffices to establish (\ref{2.26}) for $j=3$ and ends the
proof of the proposition. \hfill \qed
\bigskip

{\bf Remark:}
For $H\in]\frac{1}{2},1[$, $F(W)=(1,W^1)$ is a geometric rough path of roughness $p$,
with $pH>1$. The large deviation principle stated in Theorem \ref{t1} also holds for these
values of the parameters $H$ and $p$. Indeed, it is a consequence of (\ref{2.6}) and (\ref{2.26})
for $j=1$.


\section{Geometric rough paths on the reproducing kernel Hilbert space}

Proposition \ref{p2} implies the existence of a geometric rough path
of roughness $p$, for any $p\in]1,4[$ with $pH>1$,
lying above $h\in\h$. In this section we give
a representation of this object in terms of multiple integrals based on $h$.

We start by introducing the type of integrals to be used. They are a sort of deterministic counterpart
of the stochastic integral with respect to the fractional Brownian motion introduced in
\cite{amn} (see also \cite{ccm}).

Following \cite{amn}, for a function
$\ffi:[0,1]\to \re$ we set
\begin{align}
&||\ffi||_{K}^2=\int_0^1 \ffi(s)^2K(1,s)^2 ds +\int_0^1 ds\left(\int_s^1
|\ffi(t)-\ffi(s)|\,|K|(dt,s)\right)^2,\label{3.1}\\
&K^*(\ffi\1_{[0,t]})(s) = \ffi(s)K(t,s) + \int_s^t(\ffi(r)-\ffi(s))K(dr,s).\label{3.2}
\end{align}
Notice that,  $||\ffi||_{K}<\infty$ implies $||\ffi\1_{[0,t]}||_{K}<\infty$ as well, for any
$t\in[0,1]$. We denote by $\h_{K}$ the completion of the set $\mathcal{E}$ of
step functions on $[0,1]$ with respect to the semi-norm $||\cdot||_{K}$.

In the sequel, we set $(K\dot h)(t)=h(t)$ for any $h\in\h$ with representation given in (\ref{h})
in terms of $\dot{h}\in L^2([0,1])$.
By Lemma 1 in \cite{amn}, for any step function $\ffi\in\mathcal{E}$, we have
\beq
\label{3.3}
\int_0^1 \ffi(t) h(dt)=\int_0^1\ffi(t)(K\dot h)(dt)=\int_0^1 K^*(\ffi)(t) \dot h(t)dt.
\eeq
Thus, the linear continuous functional $\ffi\mapsto \int_0^1\ffi(t) h(dt)$
defined on $\mathcal{E}$ -endowed with the topology induced by the semi-norm $||\cdot||_{K}$-
taking values in $\re$, can be extended to $\h_K$. Hence, we attach a meaning to the indefinite
integral of $\ffi\in \h_K$ with respect to $h\in \h$ by means of the formula
\beq
\label{3.4}
\int_0^t \ffi(s) h(ds) = \int_0^1 K^*(\ffi\1_{[0,t]})(s) \; \dot h(s)ds.
\eeq
The following lemma establishes the existence of the indefinite multiple It\^o-Wiener integral with
respect to the fractional Brownian motion and its continuity. Recall that we set $K(t,s)=0$ for
$s\geq t$.

\begin{prop}
\label{p3}
Let $\ffi$ be a H\"older continuous real-valued function defined on $[0,1]$, of order
$\lambda\in(0,1)$ with $\lambda+H>\frac{1}{2}$. Then $\ffi\in\h_K$ and the function
$t\mapsto \int_0^t \ffi(s) h(ds)$ is H\"older continuous of order $H$.
\end{prop}

\noindent{\bf Proof:} First we prove that $||\ffi||_K<\infty$. Clearly,
\beqn
\int_0^1 \ffi(s)^2 K(1,s)^2 ds \le ||\ffi||_\infty^2 \int_0^1 K(1,s)^2 ds = ||\ffi||_\infty^2 <\infty,
\eeqn
with $|| \ffi||_{\infty} = \sup_{t\in [0,1]} |\ffi(t)|$. Moreover, (\ref{a.1}) implies
\beqn
\int_0^1 ds \left(\int_s^1|\ffi(t)-\ffi(s)|\,|K|(dt,s)\right)^2\le C\int_0^1 ds \left(\int_s^1
|t-s|^{\lambda+H-\frac{3}{2}}dt\right)^2
<\infty.
\eeqn
The two above inequalities yield $||\ffi||_K<\infty$.

For any $m\ge 1$, we consider the step function
$\ffi_m(s)= \sum_{l=1}^{2^m}\1_{\Del_l^m}(s)\ffi(t_{l-1}^m)$.
Since $\ffi$ is H\"older continuous,
\beqn
\sup_{l=1,\dots,2^m} \sup_{s\in\Del_l^m}|\ffi_m(s)-\ffi(s)|\le C 2^{-\lambda m}.
\eeqn
Consequently,
\beqn
\lim_{m\to\infty}\int_0^1 \left|\ffi_m(s)-\ffi(s)\right|^2 K(1,s)^2 ds
\leq \lim_{m\to\infty} C 2^{-2\lambda m}\int_0^1K(1,s)^2 ds=0.
\eeqn
Moreover,
\beqn
\lim_{m\to\infty}\big|\big(\ffi_m(t)-\ffi_m(s)\big)-\big(\ffi(t)-\ffi(s)\big)\big| \le
C\lim_{m\to\infty}2^{-\lambda m}=0,
\eeqn
and
\beqn
\sup_{m\ge 1}\big|\big(\ffi_m(t)-\ffi_m(s)\big)-\big(\ffi(t)-\ffi(s)\big)\big|\le C|t-s|^\lambda,
\eeqn
whenever $s\in\Del_l^m$, $t\in\Del_{l'}^m$, with $|l-l'|>1$.

Set
\beqn
I_m(\ffi)=\int_0^1 ds
\left(\int_s^1\big|\big(\ffi_m(t)-\ffi_m(s)\big)-\big(\ffi(t)-\ffi(s)\big)\big|\, |K|(dt,s)\right)^2
\eeqn
and $I_m(\ffi)\leq C \sum_{i=1}^3I_m^i(\ffi)$ with
\begin{align*}
I_m^1(\ffi)&=\sum_{l=1}^{2^m}\int_{\Del_l^m}ds
\left(\int_{t_{l+1}^m}^1\big|\big(\ffi_m(t)-\ffi_m(s)\big)-\big(\ffi(t)-\ffi(s)\big)\big|\,
|K|(dt,s)\right)^2,\\
I_m^2(\ffi)&=\sum_{l=1}^{2^m}\int_{\Del_l^m}ds
\left(\int_{s}^{t_l^m}\big|\ffi(t)-\ffi(s)\big|\,
|K|(dt,s)\right)^2,\\
I_m^3(\ffi)&=\sum_{l=1}^{2^m}\int_{\Del_l^m}ds
\left(\int_{\Del_{l+1}^m}\big|\big(\ffi_m(t)-\ffi_m(s)\big)-\big(\ffi(t)-\ffi(s)\big)\big|\,
|K|(dt,s)\right)^2.
\end{align*}

By the bounded convergence theorem applied first to the integral with respect to the
measure $|K|(dt,s)$ and then to the Lebesgue measure, we have
$\lim_{m\to\infty}I_m^1(\ffi)=0$.

Moreover,
\begin{align*}
I_m^2(\ffi)&\le C\sum_{l=1}^{2^m}\int_{\Del_l^m}ds \left(\int_{s}^{t_l^m}
|t-s|^{\lambda+H-\frac{3}{2}}dt\right)^2\\
&\le C 2^{-m(2\lambda+2H-1)}.
\end{align*}
Thus,
$\lim_{m\to\infty}I_m^2(\ffi)=0$.

Since
\beqn
\sup_{s,t\in[0,1]}\left(|\ffi_m(t)-\ffi(t)|+|\ffi_m(s)-\ffi(s)|\right)\le C 2^{-\lambda m},
\eeqn
it follows that
\begin{align*}
I_m^3(\ffi)&\le C\sum_{l=1}^{2^m}2^{-2\lambda m}
\int_{\Del_l^m}ds\left(\int_{\Del_{l+1}^m}|t-s|^{H-\frac{3}{2}}dt\right)^2\\
&\le C 2^{-m(2\lambda +2H-1)},
\end{align*}
and therefore, 
$\lim_{m\to\infty}I_m^3(\ffi)=0$.

Therefore, $\lim_{m\to\infty}I_m(\ffi)=0$ and we have thus established that
$\ffi\in\h_K$.

Let us now prove the H\"older continuity of the indefinite integral $\int_0^t \ffi(s) h(ds)$. Fix
$0\le t_1\le t_2\le1$. By virtue of (\ref{3.4}) and (\ref{3.2}),
\beqn
\int_0^{t_2} \ffi(s) h(ds)-\int_0^{t_1} \ffi(s) h(ds)=\sum_{i=1}^3T_i(t_1,t_2),
\eeqn
with
\begin{align*}
T_1(t_1,t_2)&=\int_0^{t_1} ds\, \dot h(s)\left(\int_{t_1}^{t_2}\ffi(r)K(dr,s)\right),\\
T_2(t_1,t_2)&=\int_{t_1}^{t_2}\ffi(s) K(t_2,s)\,\dot h(s)ds,\\
T_3(t_1,t_2)&=\int_{t_1}^{t_2} ds\, \dot h(s)\left(\int_s^{t_2}(\ffi(r)-\ffi(s))K(dr,s)\right).\\
\end{align*}
Schwarz's inequality yields
\begin{align}
|T_1(t_1,t_2)|&\le \|\ffi\|_{\infty}\, ||\dot h||_2\left(\int_0^{t_1} ds
\left(\int_{t_1}^{t_2}|K|(dr,s)\right)^2\right)^{\frac{1}{2}}\nonumber\\
&\le  \|\ffi\|_{\infty}\,
||\dot h||_2\left(\int_0^{t_1} |K(t_2,s)-K(t_1,s)|^2\, ds\right)^{\frac{1}{2}} \nonumber\\
&\leq  \|\ffi\|_{\infty}\, ||\dot h||_2\; |t_2-t_1|^H\label{3.5}.
\end{align}
Similarly,
\begin{align}
|T_2(t_1,t_2)|&\le ||\ffi||_{\infty}||\dot h||_2\left(\int_{t_1}^{t_2} K(t_2,s)^2 
ds\right)^{\frac{1}{2}}\nonumber\\
&\le ||\ffi||_{\infty}||\dot h||_2\;  |t_2-t_1|^H \label{3.6}.
\end{align}
The H\"older continuity of the function $\ffi$ together with
the upper bound given in (\ref{a.1}), imply
\begin{align}
|T_3(t_1,t_2)|&\le C||\dot h||_2 \left(\int_{t_1}^{t_2} ds
\left(\int_s^{t_2}|r-s|^{\lambda+H-\frac{3}{2}}\, dr \right)^2\right)^{\frac{1}{2}}\nonumber\\
&\le C||\dot h||_2 \; |t_2-t_1|^{\lambda+H}\label{3.7}.
\end{align}

With (\ref{3.5})--(\ref{3.7}), we have
\beqn
\left|\int_0^{t_2} \ffi(s) h(ds)-\int_0^{t_1} \ffi(s) h(ds)\right|\le C||\dot h||_2\; |t_2-t_1|^H.
\eeqn
This completes the proof of the proposition. \hfill \qed
\medskip

The preceding proposition provides a background to define indefinite iterated
integrals with respect to elements of the reproducing kernel Hilbert space of
the fractional Brownian motion, as follows.

\begin{corol}
\label{c1}
The reproducing kernel Hilbert space ${\mathcal H}$ of the fractional Brownian motion with Hurst parameter
$H\in]\frac{1}{4}, \frac{1}{2}[$ is contained in $\h_K$. Given $h\in {\mathcal H}$,
the indefinite integral $h^2_{0,t}:=\int_0^t h(s) h(ds)$
defines a $H$-H\"older continuous function. Therefore, the function $t\mapsto h^2_{0,t}$
belongs to $\h_K$. Thus, it can be integrated again
with respect to
$h$. The resulting integral inherits the $H$-H\"older continuity property.
\end{corol}

Let $g$ be a measurable Lebesgue integrable function defined on $[0,1]$.
For $l\in \{1, \, \cdots, \, 2^m\}$, set
\[a_l^m(t)= 2^m\int_{\Del_l^m\cap[0,t]} g(s)ds.\]
 Consider the linear interpolation
of $h$, that is the function $h(m)$ defined in (\ref{h(m)}). Obviously,
\begin{align}
\int_0^t g(s) h(m)(ds)& = \sum_{l=1}^{\tm +1} a_l^m(t) \Del_l^m h\nonumber\\
&=\sum_{l=1}^{\tm +1} a_l^m(t) \left((K\dot h)(t_l^m)-(K\dot h)(t_{l-1}^m)\right).\label{3.70}
\end{align}
Following the steps of the proof of Lemma 1 in \cite{amn}, consisting actually into an integration by
parts, we obtain
\beq
\label{3.8}
\int_0^t g(s) h(m)(ds)=\int_0^1K(m)^*\big(g\1_{[0,t]}\big)(s)\, \dot h(s)ds,
\eeq
with
\begin{align}
K(m)^*&\big(g\1_{[0,t]}\big)(s)=\sum_{l=1}^{\tm +1} \1_{\Del_l^m}(s) a_l^m(t)  K(t_{\tm
+1}^m,s)\nonumber\\ &+\sum_{l=1}^{\tm}\1_{\Del_l^m}(s)\sum_{l'=l+1}^{\tm
+1}(a_{l'}^m(t)-a_l^m(t))\left(K(t_{l'}^m,s)-K(t_{l'-1}^m,s)\right)\label{3.9}.
\end{align}
Notice the similarity between the expressions (\ref{3.9}) and (\ref{3.2}).
\smallskip

Our next aim is to prove that $\lim_{m\to\infty}\int_0^t G(m)(s)h(m)(ds)=\int_0^t G(s)h(ds)$, for
the pairs $G(m)=h(m)$, $G=h$, and $G(m)=\int_0^{\cdot} h(m)(s)h(m)(ds)$, $G=\int_0^{\cdot} h(s)h(ds)$,
respectively. As a consequence we shall obtain in Theorem \ref{t2} an integral expression for the
geometric rough path lying above $h$. A basic ingredient of its proof is provided by the next statement.
\begin{prop}
\label{p4}
Let $g$ be a $\lambda$-H\"older continuous real-valued function defined on $[0,1]$, with
$\lambda+H>\frac{1}{2}$. Then, there exists  a constant $C>0$ such that  for any
$0\le t_1\le t_2\le 1$, 
\beq
\label{3.10}
\sup_{m\in {\mathbb N}} \, \left|\int_0^{t_2} g(s) h(m)(ds)-\int_0^{t_1} g(s) h(m)(ds)\right|\le
C  |t_2-t_1|^H .   
\eeq
In particular, each indefinite integral $\int_0^{\cdot} g(s) h(m)(ds)$ defines a $H$-H\"older continuous
function.
\end{prop}
\noindent{\bf Proof:} Fix $m\ge 1$. Assume first that $\tmm=\tmmm$, so that
$|t_2-t_1|\leq 2^{-m}$. Owing to (\ref{3.70}),
\begin{align}
\Big|\int_0^{t_2} &  g(s) h(m)(ds)-\int_0^{t_1} g(s) h(m)(ds)\Big|\ =
\Big| 2^m\,\left( \int_{t_1}^{t_2} g(r)\, dr\right)\nonumber\\
&\qquad \times  \int_0^1 \big( K(t_{\tmm +1}^m, s) -K(t_{\tmm}^m, s)\big)\, \dot{h}(s)\, ds \Big|\nonumber \\
& \leq\|g\|_{\infty}\, \|\dot{h}\|_2\, 2^{-m(H-1)}\, |t_2-t_1|\, 
\leq  C\, |t_2-t_1|^H\, .\label{egal}
\end{align}

Suppose now that $\tmm < \tmmm$. Then using (\ref{3.70})
we have
\[ \left|\int_0^{t_2} g(s) h(m)(ds)-\int_0^{t_1} g(s) h(m)(ds)\right|\le
\sum_{j=1}^3 \left| \int_0^1 S_j(t_1,t_2,s)\, \dot{h}(s)\, ds \right|,\]
with
\begin{eqnarray*}
S_1(t_1,t_2,s)&=&\left(2^m \int_{t_1}^{t_{\tmm +1}^m} g(r) dr\right) \,
\big( K(t_{\tmm +1}^m, s) -K(t_{\tmm}^m, s)\big), \\
S_2(t_1,t_2,s)&=&\left(2^m \int_{t_{\tmmm}^m}^{t_2} g(r) dr\right) \,
\big( K(t_{\tmmm +1}^m, s) -K(t_{\tmmm}^m, s)\big), \\
S_3(t_1,t_2,s)&=& \sum_{l=\tmm +2}^{\tmmm} a^m_l(t_2)\, \Delta^m_l K(.,s)\, ,
\end{eqnarray*}
with the convention that $\sum_{l=I}^J x_l=0 $ if $I>J$.

The arguments used to prove (\ref{egal}) show that
\begin{equation} \label{j=12}
 \sum_{j=1}^2 \left| \int_0^1 S_j(t_1,t_2,s)\, \dot{h}(s)\, ds\right|  \leq C\, |t_2-t_1|^H\, .
\end{equation}
The inequalities (\ref{egal}) and  (\ref{j=12}) prove (\ref{3.10})
if $\tmmm = \tmm +1$. In order to conclude the proof, assume that $\tmmm \geq \tmm +2$ and
let us estimate $S_3(t_1,t_2,s)$.

Following again the steps of the proof of Lemma 1 in \cite{amn}, we deduce
$S_3(t_1,t_2,s) = \sum_{j=1}^3 S_{3,j}(t_1,t_2,s)$, with
\begin{align*}
S_{3,1}(t_1,t_2,s) \, =&\, 1_{]0, t^m_{\tmm +1}[}(s)\, \sum_{l=\tmm +2}^{\tmmm} a^m_l(t_2) \Delta^m_l K(.,s),\\
S_{3,2}(t_1,t_2,s)\, =&\, \sum_{l=\tmm +2}^{\tmmm } 1_{\Delta^m_l}(s)\,
a^m_l(t_2) K(t^m_{\tmmm},s), \\
S_{3,3}(t_1,t_2,s)\, =& \, \sum_{l=\tmm +2}^{\tmmm -1} 1_{\Delta^m_l}(s)\,
\Big( \sum_{l'=l+1}^{\tmmm } \big(a^m_{l'}(t_2)-a^m_l(t_2)\big)\,  \Delta^m_{l'}K(.,s)\Big)\, .
\end{align*}
Note that this decomposition is similar to that used to prove the
H\"older regularity of the indefinite stochastic integral
$\int_0^. \varphi(s)\, h(ds)$. Actually, out of the factor $\dot h$, $S_{3,j}(t_1,t_2,s)$,
$j=1,2,3$, are the analogue of the integrands of $T_j(t_1,t_2)$, $j=1,2,3$, respectively.

By (\ref{a.1}), we have that the
function $t\mapsto  K(t,s)$ is decreasing on $]s,1]$. Hence, given
$1\leq I<J$ and $s\leq t^m_{I-1}$,
\[ \sum_{l=I}^J \big|\Delta^m_lK(.,s)\big| = \big| K(t^m_J,s) - K(t^m_{I-1},s)\big|\, .\]
Since $\sup_{m,l,t} |a^m_l(t)|\leq ||g||_{\infty}$, we have
\begin{align}\label{S31}
\Big| \int_0^1 & S_{3,1}(t_1,t_2,s)\, \dot{h}(s)ds\Big|\leq  \|g\|_{\infty}\, \|\dot{h}\|_2 \nonumber \\
& \times \left( \int_0^{t^m_{\tmm +1}} ds
\big| K(t^m_{\tmmm},s) -K(t^m_{\tmm +1},s)\big|^2\right)^{\frac{1}{2}}\nonumber \\
\leq & C\left| t^m_{\tmmm}-t^m_{\tmm +1}\right|^H \leq C\, |t_2-t_1|^H 
\end{align}
and
\begin{align}\label{S32}
\Big| \int_0^1  S_{3,2}(t_1,t_2,s)\, \dot{h}(s)ds\Big| & \leq  \|g\|_{\infty}\, \|\dot{h}\|_2\nonumber 
\left( \int_{t^m_{\tmm +1}}^{t^m_{\tmmm}} ds
\big| K(t^m_{\tmmm},s)\big|^2\right)^{\frac{1}{2}}\nonumber \\
&\leq  C\left| t^m_{\tmmm}-t^m_{\tmm +1}\right|^H\leq   C\left| t_2-t_1\right|^H .
\end{align}
The H\"older continuity of $g$ implies that
for $s\in \Delta_l^m$ , $r\in \Delta_{l'}^m$ with
$[2^m t_1]+2 \leq l < l' \leq [2^m t_2]$,
$| a^m_{l'}(t_2) - a^m_l(t_2) | \leq C \left( (l'-l) 2^{-m} \right)^{\lambda}
\leq C \big(2^{-m\lambda} 1_{\{l'=l+1\}}$\linebreak $ + |r-s|^\lambda 1_{\{ l'> l+1\}} \big)$.
Therefore, since $|t_2-t_1| \geq 2^{-m}$,
\begin{align}
\Big|&\int_0^1 S_{3,3}(t_1,t_2,s)\, \dot{h}(s)ds\Big|\leq C||\dot h||_2\Big(\int_0^1\,ds\,
\sum_{l=\tmm + 2}^{\tmmm -1}\1_{\Delta_l^m}(s)\nonumber\\
&\quad \times\big(\big[\int_{t_{l+1}^m}^{t_{[2^mt_2]}^m} |r-s|^{\lambda+H-\frac{3}{2}} dr\big]^2
+\big[2^{-m\lambda}\int_{t_l^m}^{t_{l+1}^m}|r-s|^{H-\frac{3}{2}} dr\big]^2\big)\Big)^{\frac{1}{2}}\nonumber\\
&\le C ||\dot h||_2 \Big(\int_0^1\,ds\,\sum_{l=\tmm + 2}^{\tmmm -1}\1_{\Delta_l^m}(s)
\big[(t_{[2^mt_2]}^m-s)^{2\lambda+2H-1}\nonumber\\
&\quad +2^{-2 m\lambda}(t_l^m-s)^{2H-1}\big]\Big)^{\frac{1}{2}}\nonumber\\
&\le C ||\dot h||_2 \big(|t_2-t_1|^{\lambda+H}+2^{-m(\lambda+H-\frac{1}{2})}
|t_2-t_1|^{\frac{1}{2}}\big)\nonumber\\
&\le C ||\dot h||_2|t_2-t_1|^{\lambda +H}.\label{S33}
\end{align}
The inequalities (\ref{j=12})-(\ref{S33}) conclude the proof of the proposition.\hfill\qed
\smallskip

We next prove the announced result on convergence of integrals.

\begin{prop}
\label{p5} Let $h=K\dot{h}\in {\mathcal H}$, $G(m)$, $ m\ge 1$, and $G$ be real continuous
 functions defined on $[0,1]$. Assume that for any $0\le t_1\le t_2\le 1$,
\begin{description}
\item{(i)}\; $|G(t_2)-G(t_1)|\le C |t_2-t_1|^H$,
\item{(ii)}\, $ \sup_{m\in\mathbb{N}}\,|G(m)(t_2)-G(m)(t_1)|\le C\,|t_2-t_1|^H$,
\item{(iii)}\, $c(m):=\sup_{t\in[0,1]}|G(m)(t)-G(t)|\to 0$, as $m\to\infty$.
\end{description}
Then
\beqn
\label{3.11}
\lim_{m\to\infty}\sup_{t\in[0,1]} \left| \int_0^t G(m)(s) h(m)(ds) - \int_0^t G(s) h(ds)\right| =0 .
\eeqn
\end{prop}
\noindent{\bf Proof:} Set $A_l^m(t)=2^m\int_{\Del_l^m\cap[0,t]} G(m)(s) ds$. We first prove that
\beq
\label{3.12}
\lim_{m\to\infty}\sup_{t\in[0,1]}\Big|\int_0^t ds\,\dot h(s)\Big(\sum_{l=1}^{\tm +1}\1_{\Del_l^m}(s)
A_l^m(t) K(t_{\tm +1}^m,s)\Big)- G(s)K(t,s)\Big|= 0.
\eeq
Indeed, Schwarz's inequality yields
\begin{align}
&\left|\int_0^t ds\,\dot h(s)\sum_{l=1}^{\tm +1}\1_{\Del_l^m}(s)G(s)\left(K(t_{\tm
+1}^m,s)-K(t,s)\right)\right|\nonumber\\
& \leq ||\dot h||_2||G||_{\infty}\left(\int_0^1 ds \sum_{l=1}^{\tm
+1}\1_{\Del_l^m}(s)\left|K(t_{\tm
+1}^m,s)-K(t,s)\right|^2\right)^{\frac{1}{2}} \nonumber\\
& \le C \, 2^{-mH}.\label{3.13}
\end{align}
Owing to $(ii)$ and $(iii)$,
\beqn
\sup_{1\le l\le \tm}\sup_{r,s\in\Del_l^m}|G(m)(r)-G(s)|\le C(2^{-mH}+c(m)).
\eeqn
Therefore,
\begin{align}
&\left|\int_0^t ds\,\dot h(s) K(t_{\tm +1}^m,s)\left(\sum_{l=1}^{\tm}\1_{\Del_l^m}(s)
\left(A_l^m(t)-2^m\int_{\Del_l^m}G(s) dr\right)\right)\right|\nonumber\\
&\le C\big(2^{-mH}+c(m)\big) \; \int_0^1 |\dot h(s)| |K(t_{\tm +1}^m,s)|\, ds\nonumber\\
&\le C||\dot h||_2 (2^{-mH}+c(m)).\label{3.14}
\end{align}
By $(ii)$ and $(iii)$,
\beqn
\sup_{r,s\in\Del_{\tm+1}^m}\left| G(s) - 2^m  \int_{t_{\tm}^m}^t G(m)(r) dr \right| 
\le C(2^{-mH}+c(m)+||G||_{\infty})\le C.
\eeqn
Hence,
\begin{align}
&\left|\int_0^t ds\,\dot h(s) K(t_{\tm+1}^m,s)\left(\1_{\Del_{\tm +1}^m}(s)\left(2^m
\int_{\tm}^t G(m)(r) dr-G(s)\right)\right)\right|\nonumber\\
&\le C\,  \int_0^1 \1_{\Del_{\tm +1}^m}(s) |\dot h(s )| |K(t_{\tm+1}^m,s)|\,
ds
\, \le\,  C 2^{-mH}.\label{3.15}
\end{align}
With (\ref{3.13})--(\ref{3.15}), we have proved (\ref{3.12}).

The second step of the proof consists in checking that
\begin{align}
&\lim_{m\to\infty}\sup_{t\in[0,1]}\Big|\int_0^t ds\,\dot h(s)\Big(\int_s^t(G(r)-G(s))K(dr,s)\nonumber\\
&\quad-\Big[
\sum_{l=1}^{\tm}\1_{\Del_l^m}(s)\sum_{l'=l+1}^{\tm+1}\big(A_{l'}^m(t)-A_l^m(t)\big)\big(K(t_{l'}^m,s)-
K(t_{l'-1}^m,s)\Big]\Big)\Big|=0.\label{3.16}
\end{align}
Clearly,
\beqn
\int_0^t ds\,\dot h(s)\Big(\int_s^t(G(r)-G(s))K(dr,s)\Big)=\sum_{i=1}^3R_i^m(t),
\eeqn
with
\begin{align*}
R_1^m(t)&=\int_{\Del_{\tm+1}^m\cap[0,t]}ds\,\dot h(s)\Big(\int_s^t(G(r)-G(s))K(dr,s)\Big),\\
R_2^m(t)&=\sum_{l=1}^{\tm}\int_0^t ds \1_{\Delta_l^m}(s) \dot h(s)\Big(\int_s^{t_l^m}(G(r)-G(s))
K(dr,s)\Big),\\
R_3^m(t)&=\sum_{l=1}^{\tm}\int_0^t ds \1_{\Delta_l^m}(s) \dot h(s)
\Big(\sum_{l'=l+1}^{\tm+1}\int_{t_{l'-1}^m}^{t_{l'}^m}(G(r)-G(s))\1_{\{s\le r\le t\}} K(dr,s)\Big).
\end{align*}
By virtue of Schwarz's inequality,  assumption $(i)$ and (\ref{a.1}) we have
\begin{align*}
|R_1^m(t)|&\le C||\dot h||_2\Big(\int_{\Del_{\tm+1}^m\cap[0,t]}ds\Big(\int_s^t
|r-s|^{2H-\frac{3}{2}}dr\Big)^2\Big)^{\frac{1}{2}}\\
&\le C||\dot h||_2 2^{-2mH}.
\end{align*}
Using again (i) and (\ref{a.1}) we obtain
\beqn
\Big|\int_s^{t_l^m}(G(r)-G(s)) K(dr,s)\Big|\le C
(t_l^m-s)^{2H-\frac{1}{2}}.
\eeqn
It follows that
\beqn
\sup_{l=1,\dots,\tm}\sup_{s\in\Del_l^m}\Big|\int_s^{t_l^m}(G(r)-G(s))
K(dr,s)\Big|\le C 2^{-m(2H-\frac{1}{2})}
\eeqn
and consequently
\beqn 
|R_2^m(t)|\le
C2^{-m(2H-\frac{1}{2})}\sum_{l=1}^{\tm}\int_0^t ds
\1_{\Delta_l^m}(s) |\dot h(s)| \le C2^{-m(2H-\frac{1}{2})}.
\eeqn
Thus we have shown
\beqn
\lim_{m\to\infty}\sup_{t\in[0,1]}\big(|R_1^m(t)|+|R_2^m(t)|\big)=0.
\eeqn
Therefore, the proof of (\ref{3.16}) reduces to that of
\begin{align}
\lim_{m\to\infty}\sup_{t\in[0,1]} & \Big|\sum_{l=1}^{\tm}\int_0^t
ds \1_{\Delta_l^m}(s) \dot h(s)  \Big(\sum_{l'=l+1}^{\tm+1}\big(
\int_{t_{l'-1}^m}^{t_{l'}^m}\big(G(r)-G(s)\big) \1_{\{s\le r\le
t\}}\nonumber\\ &-\big(A_{l'}^m(t)-A_l^m(t)\big)\big)
K(dr,s)\Big)\Big|=0.\label{3.17}
\end{align}
Set
\begin{align*}
R_4^m(t)= & \sum_{l=1}^{\tm}\int_0^t ds \int_0^1 K(dr,s) \dot h(s) \1_{\Del_l^m}(s)
\1_{\Del_{l+1}^m}(r)\\
&\times \Big[ \big( G(r)-G(s)\big)\1_{\{s\le r\le
t\}}-\big(A_{l+1}^m(t)-A_l^m(t)\big)\Big].
\end{align*}

The H\"older continuity of $G$ together with (\ref{a.1}) implies
\begin{align}
\Big|\sum_{l=1}^{\tm} & \int_{\Del_l^m} ds \dot h(s)  \int_{\Del_{l+1}^m}  K(dr,s)
\big( G(r)-G(s)\big)\1_{\{s\le r\le t\}}\Big|\nonumber\\
&\le C \Big|\sum_{l=1}^{\tm}\int_{\Del_l^m} ds \dot
h(s)\int_{\Del_{l+1}^m}|r-s|^{2H-\frac{3}{2}} dr\Big| 
\le C\,  ||\dot h||_2 \;
2^{-m(2H-\frac{1}{2})}.\label{3.18}
\end{align}
By assumption $(ii)$, for any $l=1,\dots,\tm -1$, we have
\beqn
\sup_{t\in[0,1]}\big|A_{l+1}^m(t)-A_l^m(t)\big|\le C2^{-mH},
\eeqn
while for $l= \tm$,
\beqn
\sup_{t\in[0,1]}\big|A_{l+1}^m(t)-A_l^m(t)\big|\le C\big(2^{-mH}+c(m)+||G||_{\infty}\big)\leq C .
\eeqn
Therefore,
\begin{align}
\Big|\sum_{l=1}^{\tm}&\int_{\Del_l^m} ds \dot h(s)
\int_{\Del_{l+1}^m}   K(dr,s)
\big(A_{l+1}^m(t)-A_l^m(t)\big)\Big|\nonumber\\
&\le C 2^{-mH}\sum_{l=1}^{\tm -1}\int_{\Del_l^m} ds |\dot{h}(s)|\, \int_{\Del_{l+1}^m} dr
|r-s|^{H-\frac{3}{2}}\nonumber\\
&\quad +C\,  \int_{\Del_{\tm}^m} ds |\dot{h}(s)| \int_{\Del_{\tm +1}^m} dr
|r-s|^{H-\frac{3}{2}}\nonumber\\
&\leq C\, 2^{-mH}\, \|\dot{h}\|_2 \, \Big( \sum_{l=1}^{[2^m t] -1}
\int_{\Delta^m_l} ds |t^m_l-s|^{2H-1} \Big)^{\frac{1}{2}}
+ C\, \|\dot{h}\|_2 \, 2^{-mH} \nonumber \\
&\le C||\dot h||_2\big( 2^{-\frac{m}{2}(4H-1)} +  2^{-mH}\big) .\label{3.19}
\end{align}
From (\ref{3.18}) and (\ref{3.19}) we obtain
\beqn
\lim_{m\to\infty}\sup_{t\in[0,1]}|R_4^m(t)| = 0.
\eeqn
Consequently, it remains to prove that
\begin{align}
\lim_{m\to\infty}\sup_{t\in[0,1]} & \Big|\sum_{l=1}^{\tm
-1}\int_0^t ds  \1_{\Delta_l^m}(s)  \dot h(s)  \Big(
\sum_{l'=l+2}^{\tm+1}\Big( \int_{\Del_{l'}^m}\Big(
\big(G(r)-G(s)\big)\1_{\{s\le r\le t\}}\nonumber\\
&-\big(A_{l'}^m(t)-A_l^m(t)\big)\Big)K(dr,s)\Big)\Big|=0 .
\label{3.20}
\end{align}

Set
\beqn
\Psi_m(t;s,r)=\sum_{l=1}^{\tm-1}\sum_{l'=l+2}^{\tm}\1_{\Delta_l^m}(s)\1_{\Delta_{l'}^m}(r)
\Big( \big(G(r)-G(s)\big)-\big(A_{l'}^m(t)-A_l^m(t)\big)\Big).
\eeqn
A simple analysis based on the hypotheses $(i)-(iii)$ gives
\begin{align*}
\sup_{l=1,\dots,\tm-1} &\big|\1_{\Del_l^m}(s)\, \big(G(s)-A_l^m(t)\big)\big| 
+\sup_{l'=l+2,\dots,\tm}\big|\1_{\Delta_{l'}^m}(r)\big(G(r)-A_{l'}^m(t)\big)\big|\\
& \le C(2^{-mH}+c(m)).
\end{align*}
Thus
\beqn
|\Psi_m(t;s,r)|\le C(2^{-mH}+c(m))\1_{\{(s,r)\in[0,t]^2: s\le r\}}
\eeqn
and therefore,
\beqn
\lim_{m\to\infty}\sup_{t\in[0,1]}\sup_{(s,r)\in[0,1]^2: s\le r}|\Psi_m(t;s,r)|=0.
\eeqn
\label{3.21}
For any $l=1,\dots,\tm -1$, $l'=l+2,\dots,\tm$,
\beqn
\1_{\Del_l^m}(s) \1_{\Delta_{l'}^m}(r)|A_{l'}^m(t)-A_l^m(t)|\le C |r-s|^H.
\eeqn
Indeed, a change of variables and the assumption $(ii)$ yield
\begin{align*}
\1_{\Del_l^m}(s) & \1_{\Delta_{l'}^m}(r)  |A_{l'}^m(t)-A_l^m(t)|=\1_{\Del_l^m}(s) \1_{\Delta_{l'}^m}(r)\\
&\quad \times 2^m \big(\big|\int_{\Delta_{l'}^m} G(m)(u) du - \int_{\Delta_{l}^m} G(m)(u) du\big|\big)\\
&\le C\, \Big(\frac{l'-l}{2^m}\Big)^H 
\le C |r-s|^H.
\end{align*}
Hence
\beqn
\sup_{m\ge 1}|\Psi_m(t;s,r)|\le C|r-s|^H\1_{\{(s,r)\in[0,t]^2: s\le r\}},
\eeqn
and consequently,
\beqn
\sup_{m\ge 1}\sup_{t\in[0,1]}|\Psi_m(t;s,r)|\le C|r-s|^H\1_{\{(s,r)\in[0,1]^2: s\le r\}}.
\eeqn
The function $(s,r)\mapsto \dot h(s)|r-s|^H\1_{\{(s,r)\in[0,1]^2: s\le r\}}$ is integrable on the set
$[0,1]^2$ with respect to the measure $ds\times |K|(dr,s)$. Hence,
\beq
\label{3.22}
\lim_{m\to\infty}\sup_{t\in[0,t]}\int_0^1 ds\dot h(s) \int_0^1 K(dr,s) \Psi_m(t;s,r)=0.
\eeq
In order to complete the proof of (\ref{3.20}), we must check that
\beq
\label{3.23}
\lim_{m\to\infty}\sup_{t\in[0,1]}|R_5^m(t)| = 0,
\eeq
where
\begin{align*}
R_5^m(t)= & \sum_{l=1}^{\tm -1}  \int_0^t ds\,  \dot h(s)\1_{\Del_l^m}(s) \int_0^1K(dr,s)\,
\1_{\Del_{\tm +1}^m}(r)\\
&\quad\times\Big( \big(G(r)-G(s)\big)\1_{\{s\le r\le t\}}-\big(A_{\tm +1}^m(t)-A_l^m(t)\big)\Big).
\end{align*}
For $l=1\dots,\tm -1$,
\begin{align*}
\1_{\Del_l^m}(s)\big|G(s)-A_l^m(t)\big|& =\1_{\Del_l^m}(s) 2^m\Big|\int_{\Del_l^m}\big(G(m)(u)-G(s)\big)
du\Big|\\
&\le C\big(2^{-mH}+c(m)\big) \leq C
\end{align*}
and
\beqn
\1_{\Del_l^m}(s)\1_{\Del_{\tm +1}^m}(r)\big|G(r)\1_{\{s\le r\le t\}}-A_{\tm +1}^m(t)\big|
\le C.
\eeqn
Moreover, given any $a\in ]0,H[$,
\begin{align*}
&\int_0^{t_{\tm-1}^m} ds \Big(\int_{\Delta^m_{[2^mt]+1}} |r-s|^{H-\frac{3}{2}} \, dr \Big)^2\\
&\le C\int_0^{t_{\tm-1}^m} ds
\big(t_{\tm}^m-s\big)^{2(1-a)(H-\frac{1}{2})}\,
\big((t^m_{[2^mt]}-s)^{H-\frac{1}{2}}-
 (t_{\tm+1}^m-s)^{H-\frac{1}{2}}\big)^{2a}\\
 &\le C 2^{-2am}\int_0^{t_{\tm-1}^m} ds (t_{\tm}^m-s)^{2H-1-2a} \leq  C 2^{-2am}.\\
\end{align*}
Hence for $a\in ]0,H[$,
\[R^m_5(t)\leq
C\, ||\dot h||_2\, 
2^{-am}. \]
This clearly implies (\ref{3.23}) and concludes the proof of the proposition. \hfill\qed

The following theorem gives an integral representation of the geometric
rough path $(1, h^1, h^2, h^3)$
associated with $h\in {\mathcal H}$.

\begin{teor}
\label{t2}
Let $h=K \dot{h}$ be an element of the reproducing kernel Hilbert space of the fractional
Brownian motion with Hurst parameter $H\in ]\frac{1}{4}, \frac{1}{2}[$. Then for every $s<t$, $i,j, \kappa
\in \{1, \cdots, d\}$,
\begin{align}
h^{2,i,j}_{s,t}=& \int_0^1 K^*\left(h^{1,i}_{0,.} \1_{]s,t]}\right)(u)\, \dot{h}^j(u)\, du\, 
-h_{0,s}^{1,i} h_{s,t}^{1,j},\label{3.24}\\
h^{3,i,j,\kappa}_{s,t}= & \int_0^1  K^*\left(h^{2,i,j}_{0,.} \1_{]s,t]}\right)(u)\,
\dot{h}^{\kappa}(u)\, du 
-h_{0,s}^{2,i,j} h_{s,t}^{1,\kappa}-h_{0,s}^{1,i}h_{s,t}^{2,j,\kappa}.\label{3.25}
\end{align}
\end{teor}
\noindent{\bf Proof:} For simplicity, we shall assume $d=1$ and consequently, we remove the indices
$i,j,\kappa$.

To prove (\ref{3.24}), set $G:=h$ and $G(m):= h(m)$. Then owing
to (\ref{2.25.pre}), $G$ is $H$-H\"older continuous. Suppose
$[2^ms]=[2^m t]$. Then $|t-s|\leq 2^{-m}$ and we have
\beqn
|h(m)(t)-h(m)(s)|\leq C\,2^{-m(1-H)} |t-s|\leq C\, |t-s|^H.
\eeqn

Assume $[2^m s]< [2^m t]$. Since 
\begin{align*}
|h(m)(t)-h(m)(s)|\leq & \,  \big|h(m)(t)-h(m)(t^m_{[2^m t]})\big| +
\big|h(m)(t^m_{[2^ms]+1})-h(m)(s)\big|\\
&+\big| h(t^m_{[2^m t]}) - h(t^m_{[2^m s]+1})\big|, 
\end{align*}
the H\"older continuity of $h$ yields
\beqn
\sup_{m\in \mathbb{N}}|h(m)(t)-h(m)(s)|\leq C \, |t-s|^H.
\eeqn  
Moreover,  $\sup_{t\in[0,1]} \big| h(m)(t) - h(t)|\leq C\, 2^{-mH}$. Thus, the
assumptions (i)-(iii) of Proposition \ref{p5} are satisfied, so
that for every $r\in [0,1]$, the sequence $(\int_0^r h(m)(u) \,h(m)(du), m\ge 1)$, 
converges to
\beqn
\int_0^r h(u) h(du) = \int_0^1 K^*\big(h\1_{]0,r]}\big)(u)\, \dot{h}(u)\,
du.
\eeqn
The construction of the geometric rough path based on $h$ given in
Proposition \ref{p2} shows that 
$h^2_{0,r}=\int_0^r h(u) h(du)$. Then, formula (\ref{3.24}) follows from the 
multiplicative properties of rough paths.

For the proof of (\ref{3.25}), we fix
$G(.):=h^2_{0,.}$ and $G(m)(.):=h(m)_{0,.}^2$.
Corollary \ref{c1}, Proposition \ref{p4} and the results set up in the first part of
this proof show that the assumptions of Proposition \ref{p5} hold true.
Therefore, for any fixed $r\in[0,1]$ the sequence $(\int_0^r h(m)^2_{0,.}(u)h(m)(du), m\ge 1)$
converges to $\int_0^r h^2_{0,.}(u) h(du)=\int_0^1 K^*(h_{0,.}^2\1_{]0,r]})(u) \dot h(u) du$.
By Proposition \ref{p2}, the limit must coincide with $h_{0,r}^3$. 
Then the expression (\ref{3.25}) follows from the
multiplicative property of rough paths.
\hfill\qed

\smallskip

\noindent {\bf Acknowledgments} This paper has been written while the first named
author was visiting the Centre de Recerca Matem\`atica in Bellaterra.

\end{document}